\RequirePackage{fix-cm}

\documentclass{article}

\usepackage{graphicx}
\usepackage{amssymb}
\usepackage{amsmath}
\usepackage{eucal}

\usepackage[driverfallback=dvipdfm,breaklinks=true,pdfcreator={LaTeX e Dvips}]{hyperref}
\usepackage{float}

\usepackage{mathptmx}
\usepackage[numbers]{natbib}
\usepackage{color}
\usepackage{ulem}
\usepackage{booktabs}
\usepackage{multirow}
\usepackage{textcomp}
\usepackage{bm}

\makeatletter
\let\MYcaption\@makecaption
\makeatother
\usepackage{subcaption}
      \makeatletter
\let\@makecaption\MYcaption
\makeatother

\title{Adjoint lattice kinetic scheme for topology optimization\\ in fluid problems}

\author{Yuta Tanabe$^\text{a,}\footnote{Corresponding author: {\tt 4524702@ed.tus.ac.jp} (Yuta Tanabe)}$,
        Kentaro Yaji$^\text{b}$ ,
        Kuniharu Ushijima$^\text{a}$ \\[12pt]
$^\text{a}$\textit{Department of Mechanical Engineering, Tokyo University of Science,}\\
              \textit{6-3-1, Niijuku,
              Katsushika-ku, Tokyo 125-8585, Japan}\\
$^\text{b}$\textit{Department of Mechanical Engineering, Osaka University,}\\
              \textit{2-1, Yamadaoka,
              Suita, Osaka 565-0871, Japan}}

\begin{document}

\maketitle

\begin{abstract}
    This paper proposes a topology optimization method for non-thermal and thermal fluid problems using the Lattice Kinetic Scheme (LKS).
    LKS, which is derived from the Lattice Boltzmann Method (LBM), requires only macroscopic values, such as fluid velocity and pressure, whereas LBM requires velocity distribution functions, thereby reducing memory requirements.
    The proposed method computes design sensitivities based on the adjoint variable method, and the adjoint equation is solved in the same manner as LKS; thus, we refer to it as the \textit{Adjoint Lattice Kinetic Scheme} (ALKS).
    A key contribution of this method is the proposed approximate treatment of boundary conditions for the adjoint equation, which is challenging to apply directly due to the characteristics of LKS boundary conditions.
    We demonstrate numerical examples for steady and unsteady problems involving non-thermal and thermal fluids, and the results are physically meaningful and consistent with previous research, exhibiting similar trends in parameter dependencies, such as the Reynolds number.
    Furthermore, the proposed method reduces memory usage by up to 75\% compared to the conventional LBM in an unsteady thermal fluid problem.
\end{abstract}

\section{Introduction}
\label{sec1}

Topology optimization is one of the structural optimization methods, and unlike size optimization or shape optimization which were previously proposed, the method can generate weight-reducing voids inside a structure.
Consequently, it has attracted significant attention due to its potential to generate structures that are not based on human experience.
Bends{\o}e and Kikuchi~\cite{bendsoe1988generating} were the first to introduce topology optimization into the structural mechanics field.
Subsequently, Borravall and Petersson~\cite{borrvall2003topology} pioneered the methodology in the fluid mechanics field by formulating it for a Stokes flow problem.
Later, it was applied to steady-state Navier-Stokes laminar flow~\cite{gersborg2005topology}\cite{olesen2006}, turbulent flow~\cite{kontoleontos2013adjoint}\cite{DILGEN2018363}, unsteady flow~\cite{kreissl2011topology}\cite{DENG20116688}, forced convection~\cite{yoon2010topological}\cite{Matsumori2013} and natural convection problems~\cite{alexandersen2014topology}.
Additionally, in more practical scenarios, this methodology was extended to three-dimensional heatsink design problems under natural convection~\cite{ALEXANDERSEN2016876}, and the designed heatsink was successfully manufactured~\cite{LEI2018396}.

In the topology optimization process, it is often necessary to solve state equations many times, and the implicit scheme based on Finite Element Method (FEM) is commonly employed for this purpose.
However, this approach requires solving large-scale simultaneous equations, which significantly increases computational costs.
One potential solution to this challenge is to use Lattice Boltzman Method (LBM)~\cite{chen1998lattice}, where fluid is modeled as a set of fictitious particles with finite velocity.
The collisions and propagations of the particles are calculated using distribution functions, and macroscopic values are computed from their moments.
Since LBM is a fully explicit scheme, it does not require solving large-scale simultaneous equations and is well-suited for parallel computing, which is a key trend in accelerating computations, as most calculations are performed independently at each lattice point.
Pingen et al.~\cite{pingen2007topology} proposed a topology optimization method that solves the state field using LBM.
Nevertheless, their approach still required solving large-scale simultaneous equations to compute design sensitivities.
Subsequently, Yaji et al.~\cite{yaji2014topology} introduced a method that calculates design sensitivities at computational costs nearly equivalent to LBM by employing adjoint LBM~\cite{krause2013adjoint}, in which adjoint equations are discretized and solved using the same approach as LBM.
This LBM-based topology optimization method has since been applied to various fields, including thermal fluid problems~\cite{YAJI2016355}\cite{DUGAST2018376}\cite{Luo2023}, unsteady problems~\cite{chen2017local}\cite{Yaji2018}\cite{NGUYEN202082}, stabilized with the MRT model~\cite{LUO2024121732} and the porous model~\cite{LUO2023106007}.

We also applied the topology optimization method using LBM to natural convection problems~\cite{Tanabe2023}, during this process, memory usage became a critical issue, especially for unsteady problems.
This is because, for unsteady problems, state field values for all time steps must be stored to compute design sensitivities.
Additionally, in LBM, these values represent velocity distribution functions, which significantly increase memory requirements.
It should be noted that this substantial memory consumption is not only observed in natural convection but also in forced convection.

Several approaches have been proposed to reduce memory usage, including the checkpoint algorithm~\cite{Griewank2000}, which achieves this by recalculating state field values, the local-in-time method~\cite{YAMALEEV20105394}, which splits the time interval and approximates design sensitivities, and a hybrid approach~\cite{Theulings2024} that combines these two methods.
In addition, the ``poor-man's'' approach~\cite{Asmussen2019}, which simplifies the governing equations of thermal fluid with certain approximations, and a method that approximates specific three-dimensional problems to two-dimensional ones~\cite{ZENG2020119681} were also proposed.

A noteworthy development in this context is the Lattice Kinetic Scheme (LKS), proposed by Inamuro et al.~\cite{Coveney2002lattice}, which reduces memory usage by computing exclusively with macroscopic quantities, such as pressure and fluid velocity, thereby eliminating the necessity of storing velocity distribution functions.
This is accomplished by tuning LBM parameters and adding specific terms.
In this study, we addressed the memory usage issue in topology optimization by employing LKS.
A topology optimization method using LKS has already been proposed by Xie et al.~\cite{XIE2021251}; however, their method closely resembles LBM for the following two reasons.
First, in their method, adjoint equations are derived from the fully continuous Boltzman equation.
During this process, the fluid velocity derivative terms in the local equilibrium distribution function, characteristic of LKS, are canceled out.
As a result, the adjoint field value derivative terms corresponding to them are not explicitly present in the adjoint equations.
Second, in their method, the form of the derived adjoint equations corresponds to the local equilibrium distribution function of LBM or LKS, and the adjoint equations are not described using only macroscopic values, as in LKS.
In contrast, in the proposed adjoint equations in our study are derived from the discrete velocity Boltzman equations and they are discretized in the same manner as in LKS, consequently form of adjoint equations are described corresponding to LKS.
The adjoint equations are solved in the same manner as the LKS, therefore, we refer to this method as the \textit{Adjoint Lattice Kinetic Scheme} (ALKS).

This paper is organized as follows:
Section~\ref{sec2} introduces the general theory of topology optimization, as well as the LBM and LKS. 
The section also includes the derivation of the ALKS and validation through comparison with finite difference approximations, followed by a brief explanation of the optimization procedure.
Section~\ref{sec3} presents numerical examples for both thermal and non-thermal problems, addressing both steady and unsteady scenarios.
In these examples, we also compare the memory consumption of the proposed method, which uses the LKS approach, with that of the conventional method, based on the LBM.
Finally, Section~\ref{sec4} concludes the paper by summarizing the key findings. 

\section{Formulation}
\label{sec2}

\subsection{Topology optimization}
\label{sec21}

In topology optimization, structural optimization is regarded as the problem of defining the boundaries of the design domein $\Omega$ and the fixed design domain $D\subset\mathbb{R}^d$ ($d$: spatial dimension), where $\Omega \subseteq D$. 
Then, the original optimization problem is replaced by a material distribution optimization problem by introducing the characteristic function $\chi_{\Omega}:D\rightarrow\{0,1\}$.
The characteristic function is defined as: 
\begin{equation}
    \chi_{\Omega}({\bm x}) =\left\{
    \begin{array}{cl}
        1 & \mbox{if ${\bm x} \in \Omega$}              \\
        0 & \mbox{if ${\bm x} \in D\backslash \Omega$}, \\
    \end{array}
    \right.
\end{equation}
where ${\bm x}$ expresses the position in $D$.
The optimization problem directly discretizing $\chi_{\Omega}$, known as 0-1 optimization, is challenging; thus, it is necessary to relax or regularize the problem.
One of the most well-known methods is the density method~\citep{bendsoe2003topology}, and its basic idea is to replace $\chi_{\Omega}$ with a continuous function, i.e., $\gamma:D\rightarrow[0,1]$, where $\gamma$ represents the design variable space.
Using $\gamma$, the topology optimization problem based on the density method is expressed as:
\begin{align}
    \begin{array}{ll}
        \underset{\gamma}{\text{minimize }}\  & J(\gamma, \bm{U}(\gamma))                                  \\
        \text{subject to }\                   & G(\gamma, \bm{U}(\gamma))\leqslant 0,
        \vspace{1.5mm}                                                                                     \\
                                              & 0\leqslant\gamma(\bm{x})\leqslant 1,\ \forall\bm{x}\in D,
    \end{array}
    \label{eq:obj}
\end{align}
where $J$ and $G$ are an objective functional and a constraint functional, respectively.
$\bm{U}$ expresses the state variables, which are solutions of partial differential equations governing physical phenomena.
In this study, $\bm{U}$ refers to the solution of steady or unsteady non-thermal or thermal fluid systems.
The governing equations are discritized using LKS, whose details are described in Section~\ref{sec24}, and the optimization problem in Eq.~\eqref{eq:obj} is solved using a gradient-based method.

In this study, the analysis domain is described as $\mathcal{O}=D\cup\Omega_\text{non}$ and $D\cap\Omega_\text{non}=\emptyset$, also $\Omega_\text{non}$ expresses non-design domain.
The analysis domain consists of fluid region $\Omega_\text{f}$ and solid region $\Omega_\text{s}$, and in density method, design variable means $\Omega_\text{s}$ when $\gamma=0$ and $\Omega_\text{f}$ when $\gamma=1$, and grayscale region, which is $0<\gamma<1$, corresponds to porous media.
State fields are solved within the time interval $\mathcal{I}=[t_0,t_1]$, in which $t_0$ and $t_1$ mean initial and final time, respectively.

\subsection{Topology optimization for thermal fluid problems}
\label{sec22}

For thermal fluid problems, state fields $\bm{U}$ coresponds to pressure $p:\Omega_\text{f}\times\mathcal{I}\rightarrow\mathbb{R}$, fluid velocity $\bm{u}:\Omega_\text{f}\times\mathcal{I}\rightarrow\mathbb{R}^d$ and temperature $T:\mathcal{O}\times\mathcal{I}\rightarrow\mathbb{R}$, and they are solutions of equation of continuity (Eq.~\eqref{eq:cont}), Navier-Stokes equation (Eq.~\eqref{eq:ns}) and energy equation (Eq.~\eqref{eq:energy}).
\begin{align}
     & \frac{\partial u_\alpha}{\partial x_\alpha}=0,\label{eq:cont}                                                                                                                                                                 \\
     & \frac{\partial u_\alpha}{\partial t}+u_\beta\frac{\partial u_\alpha}{\partial x_\beta}=-\frac{\partial p}{\partial x_\alpha}+\nu\frac{\partial^2u_\alpha}{\partial {x_\beta}^2}-\alpha_\gamma u_\alpha+G_\alpha,\label{eq:ns} \\
     & \frac{\partial T}{\partial t}+u_\alpha\frac{\partial T}{\partial x_\alpha}=\frac{\partial}{\partial x_\alpha}\left(K_\gamma\frac{\partial T}{\partial x_\alpha}\right)+Q,\label{eq:energy}
\end{align}
where $\nu$ represents kinematic viscosity, and the subscripts $\alpha$ and $\beta$ denote the $x$ and $y$ directions in two-dimensional problems, and $x$, $y$ and $z$ in three-dimensional problems.
The third term of the right-hand side in Eq.~\eqref{eq:ns} is external force based on Brinkman model, and $\alpha_\gamma$ is a function of the design variable $\gamma$, defined as:
\begin{equation}
    \alpha_\gamma=\bar{\alpha}\frac{q_\alpha\left(1-\gamma\right)}{q_\alpha+\gamma},\label{eq:brinkman}
\end{equation}
where $\bar{\alpha}$ is a large enough positive value and $q_\alpha>0$ is a parameter tuning convexity of $\alpha_\gamma$, namely, it effects only in solid region and make $u_\alpha$ being 0 asymptotically.
Moreover, the thermal diffusivity $K_\gamma$ is also a function of the design variable $\gamma$, and is expressed as:
\begin{equation}
    K_\gamma=K_\text{f}+\left(K_\text{s}-K_\text{f}\right)\frac{q_K\left(1-\gamma\right)}{q_K+\gamma},\label{eq:conductivity}
\end{equation}
where $K_\text{f}$ and $K_\text{s}$ are thermal diffusivity of fluid and solid respectively.

\subsection{LBM}
\label{sec23}

First, we explain the LBM, which serves as the basis for LKS.
In LBM, a thermal fluid is approximated as an aggregation of fictitious particles with a finite set of velocities, and its behavior is simulated through the calculation of the motion of these particles.
The finite set of velocities differs depending on the lattice gas model; in this explanation, we use ${\bm c}_0, {\bm c}_1, \cdots, {\bm c}_{Q-1}\in\mathbb{R}^d$.
The pressure, fluid velocity, and temperature fields are obtained from the moments of the velocity distribution functions $f_i:\mathcal{O}\times\mathcal{I}\rightarrow\mathbb{R}$($i=0,1,\cdots,Q-1$) and $g_i:\mathcal{O}\times\mathcal{I}\rightarrow\mathbb{R}$($i=0,1,\cdots,Q-1$) as follows:
\begin{align}
     & p=\frac{\rho}{3}=\frac{1}{3}\sum_{i=0}^{Q-1}f_i,\label{eq:lbmp} \\
     & {\bm u}=\sum_{i=0}^{Q-1}{\bm c}_if_i,\label{eq:lbmu}            \\
     & T=\sum_{i=0}^{Q-1}g_i,\label{eq:lbmT}
\end{align}
where $\rho$ is the fluid density.
The velocity distribution functions $f_i$ and $g_i$ are obtained by solving the lattice Boltzmann equations, as expressed in Eqs.~\eqref{eq:lbef}--\eqref{eq:lbeg_ex}.
\begin{align}
     & f_i^*\left({\bm x}+{\bm c}_i\Delta x,t+\Delta t\right)= f_i\left({\bm x},t\right)-\frac{1}{\tau_\text{f}}\left\{f_i\left({\bm x},t\right)-f^\text{eq}_i\left({\bm x},t\right)\right\}, \label{eq:lbef}                                                   \\
     & f_i\left(\bm{x},t+\Delta t\right)=f_i^*\left(\bm{x},t+\Delta t\right)-3\Delta xw_ic_{i\alpha}\alpha_\gamma\left({\bm x}\right)u_\alpha\left({\bm x},t+\Delta t\right)+3\Delta xw_ic_{i\alpha}G_\alpha\left({\bm x},t+\Delta t\right), \label{eq:lbef_ex} \\
     & g_i^*\left({\bm x}+{\bm c}_i\Delta x,t+\Delta t\right)= g_i\left({\bm x},t\right)-\frac{1}{\tau_\text{g}}\left\{g_i\left({\bm x},t\right)-g^\text{eq}_i\left({\bm x},t\right)\right\}, \label{eq:lbeg}                                                   \\
     & g_i\left(\bm{x},t+\Delta t\right)=g_i^*\left(\bm{x},t+\Delta t\right)+\Delta xQ\left({\bm x},t+\Delta t\right), \label{eq:lbeg_ex}
\end{align}
where $\tau_\text{f}$ and $\tau_\text{g}$ are the dimensionless relaxation time, respectively, and $f^\text{eq}_i$ and $g^\text{eq}_i$ are the local equilibrium distribution functions and expressed as:
\begin{align}
     & f^\text{eq}_i=w_i\left(\rho+3c_{i\alpha}u_\alpha+\frac{9}{2}c_{i\alpha}u_\alpha c_{i\beta}u_\beta-\frac{3}{2}u_\alpha^2\right), \label{eq:feq} \\
     & g^\text{eq}_i=w_iT\left(1+3c_{i\alpha} u_\alpha\right). \label{eq:geq}
\end{align}
$w_i$ are the weight coefficient differing depending on the lattice gas model and for example in D2Q9 model they are as 
\begin{equation}
    w_i=
    \left\{
    \begin{alignedat}{3}
        &4/9&& \left(i=0\right)\\
        &1/9&& \left(i=1,2,3,4\right)\\
        &1/36&&\left(i=5,6,7,8\right).
    \end{alignedat}
    \right. 
    \label{eq:D2Q9w}
\end{equation}

\subsection{LKS}
\label{sec24}

Next, we explain the procedure for deriving the basic equations of LKS.
By setting $\tau_\text{f}=\tau_\text{g}=1$ in Eq.~\eqref{eq:lbef} and Eq.~\eqref{eq:lbeg}, and replacing $\bm{x}$ with $\bm{x}-\bm{c}_i\Delta x$, we derive Eqs.~\eqref{eq:lbef2} and \eqref{eq:lbeg2}.
\begin{align}
     & f_i\left(\bm{x},t+\Delta t\right)=f_i^\text{eq}\left(\bm{x}-\Delta x,t\right),\label{eq:lbef2} \\
     & g_i\left(\bm{x},t+\Delta t\right)=g_i^\text{eq}\left(\bm{x}-\Delta x,t\right).\label{eq:lbeg2}
\end{align}
Therefore, the state fields are obtained by Eqs.~\eqref{eq:rho}, \eqref{eq:u}, and \eqref{eq:T}.
\begin{align}
    \rho\left(\bm{x},t\right)   & =\sum_{i=0}^{Q-1}f_i^{eq}\left(\bm{x}-\bm{c}_i\Delta x,t-\Delta t\right), \label{eq:rho}       \\
    \bm{u}\left(\bm{x},t\right) & =\sum_{i=0}^{Q-1}\bm{c}_if_i^{eq}\left(\bm{x}-\bm{c}_i\Delta x,t-\Delta t\right), \label{eq:u} \\
    T\left(\bm{x},t\right)      & =\sum_{i=0}^{Q-1}g_i^{eq}\left(\bm{x}-\bm{c}_i\Delta x,t-\Delta t\right). \label{eq:T}
\end{align}
Thus, in LKS, calculations can be performed without the need for velocity distribution functions $f_i$ or $g_i$.
In this study, we use the local equilibrium distribution function $f_i^\text{eq}$, as employed by Xie et al.~\cite{XIE2021251} based on the work of Inamuro et al.~\cite{Coveney2002lattice}, and $g_i^\text{eq}$, proposed by Inamuro et al.~\cite{Coveney2002lattice}.
\begin{align}
     & f_i^\text{eq}=w_i\left\{\rho+3c_{i\alpha}u_\alpha+\frac{9}{2}c_{i\alpha}u_\alpha c_{i\beta}u_\beta-\frac{3}{2}u_\alpha^2+\Delta xA\left(\frac{\partial u_\alpha}{\partial x_\beta}+\frac{\partial u_\beta}{\partial x_\alpha}\right)c_{i\alpha}c_{i\beta}\right\},\label{eq:feq2} \\
     & g_i^\text{eq}=w_iT\left(1+3c_{i\alpha}u_\alpha+\Delta xB_\gamma\frac{\partial T}{\partial x_\alpha}c_{i\alpha}\right).\label{eq:geq2}
\end{align}
Moreover, the relationship between $A$ and kinematic viscosity $\nu$, and between $B_\gamma$ and thermal diffusivity $K_{\gamma}$, are described by Eqs.~\eqref{eq:A} and \eqref{eq:B}, respectively.
\begin{align}
    \nu      & =\left(\frac{1}{6}-\frac{2}{9}A\right)\Delta x, \label{eq:A}        \\
    K_\gamma & =\left(\frac{1}{6}-\frac{1}{3}B_\gamma\right)\Delta x. \label{eq:B}
\end{align}

In the LKS, calculations are based solely on macroscopic values, enabling direct prescription of pressure or flow velocity at boundary conditions.
In this study, the inlet boundary condition is specified as $u_\alpha=\bar{u}_\alpha$, where $\bar{u}_\alpha$ represents a prescribed flow velocity component.
For wall boundaries, this is specifically defined as $\bar{u}_\alpha=0$.
On the other hand, the outlet boundary condition is given by $p =\bar{p}$ and $u_\text{s} =\bar{u}_\text{s}$, where $u_\text{s}$ denotes the flow velocity component parallel to the wall surface, and $\bar{p}$ and $\bar{u}_\text{s}$ are prescribed values for these boundary conditions.
Additionally, for thermal boundary conditions, a prescribed temperature boundary condition is defined as $T=\bar{T}$, while a prescribed heat flux boundary condition is expressed as $-q_\alpha n_\alpha=\bar{q}_\text{w}$, where $\bar{T}$ represents the specified temperature, $n_\alpha$ is the outward normal vector on the boundary, and $\bar{q}_\text{w}$ denotes the heat flux applied at the boundary.

\subsection{ALKS and sensitivity analysis}
\label{sec25}
In this section, we explain the \textit{Adjoint Lattice Kinetic Scheme} (ALKS).
We define the functional $J$ as follows: 
\begin{equation}
    J=\int_{\mathcal{I}}\int_\mathcal{O}J_\mathcal{O}d\Omega dt+\int_{\mathcal{I}}\int_{\partial\mathcal{O}}J_{\partial\mathcal{O}}d\Gamma dt. \label{eq:fj}
\end{equation}
Using the adjoint variable method, the sensitivity of $J$ is given by
\begin{align}
    \langle J^\prime,\delta\gamma\rangle= & \int_{\mathcal{I}}\int_\mathcal{O}\left(\frac{\partial J_\mathcal{O}}{\partial\gamma}+3\Delta x\frac{\partial\alpha_\gamma}{\partial\gamma}u_\alpha\tilde{u}_\alpha-3\Delta x\frac{\partial G_\alpha}{\partial\gamma}\tilde{u}_\alpha-\frac{1}{\varepsilon_\text{g}}\Delta x\frac{\partial B_\gamma}{\partial\gamma}\frac{\partial T}{\partial x_\alpha}\tilde{q}_\alpha-\Delta x\frac{\partial Q}{\partial\gamma}\tilde{T}\right)\delta\gamma d\Omega dt \notag \\
                                          & +\int_{\mathcal{I}}\int_{\partial\mathcal{O}}\frac{\partial J_{\partial\mathcal{O}}}{\partial\gamma}\delta\gamma d\Gamma dt,
\end{align}
where $\tilde{u}_\alpha$, $\tilde{T}$ and $\tilde{q}_\alpha$ are determined by solving the adjoint equations in a manner similar to LKS, as described below: 
\begin{align}
    \tilde{\rho}\left(\bm{x},t\right)=            & \sum_{i=0}^{Q-1}w_i\tilde{f}_i^\text{eq}\left(\bm{x}+\bm{c}_i\Delta x,t+\Delta t\right),\label{eq:ap}                      \\
    \tilde{u}_\alpha\left(\bm{x},t\right)=        & \sum_{i=0}^{Q-1}w_ic_{i\alpha}\tilde{f}_i^\text{eq}\left(\bm{x}+\bm{c}_i\Delta x,t+\Delta t\right),\label{eq:au}           \\
    \tilde{s}_{\alpha\beta}\left(\bm{x},t\right)= & \sum_{i=0}^{Q-1}w_ic_{i\alpha}c_{i\beta}\tilde{f}_i^\text{eq}\left(\bm{x}+\bm{c}_i\Delta x,t+\Delta t\right),\label{eq:as} \\
    \tilde{T}\left(\bm{x},t\right)=               & \sum_{i=0}^{Q-1}w_i\tilde{g}_i^\text{eq}\left(\bm{x}+\bm{c}_i\Delta x,t+\Delta t\right),\label{eq:aT}                      \\
    \tilde{q}_\alpha\left(\bm{x},t\right)=        & \sum_{i=0}^{Q-1}w_ic_{i\alpha}\tilde{g}_i^\text{eq}\left(\bm{x}+\bm{c}_i\Delta x,t+\Delta t\right). \label{eq:aq}
\end{align}
Here, $\tilde{f}_i^{eq}$ and $\tilde{g}_i^{eq}$ are expressed in Eqs.~\eqref{eq:afeq} and \eqref{eq:ageq}, respectively.
\begin{align}
    \tilde{f}_i^{eq} & =\tilde{\rho}+3c_{i\alpha}\left(\tilde{u}_\alpha+3\tilde{s}_{\alpha\beta}u_\beta-\tilde{\rho}u_\alpha\right)-\Delta xA\frac{\partial}{\partial x_\beta}\left(\tilde{s}_{\alpha\beta}+\tilde{s}_{\beta\alpha}\right)c_{i\alpha}, \label{eq:afeq} \\
    \tilde{g}_i^{eq} & =\tilde{T}+3\tilde{q}_\alpha u_\alpha-\Delta xB_\gamma\frac{\partial \tilde{q}_\alpha}{\partial x_\alpha}-\Delta x\frac{\partial B_\gamma}{\partial x_\alpha}\tilde{q}_\alpha. \label{eq:ageq}
\end{align}

\begin{figure}[t]
    \centering
    \includegraphics[width=100mm]{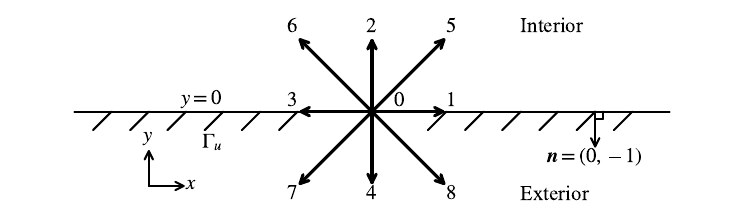}
    \caption{Boundary condition at $y=0$ using D2Q9 model}
    \label{fig:boundary_condition}
\end{figure}

Next, we describe the boundary conditions for ALKS.
Although the details of the derivation are in appendix A, on the inlet boundary, Eq.~\eqref{eq:afboundaryu} must be satisfied.
\begin{equation}
    n_\alpha\left\{\tilde{f}_i\delta_{\alpha\beta}-\frac{1}{\varepsilon_\text{f}}\Delta xA\left(\tilde{s}_{\alpha\beta}+\tilde{s}_{\beta\alpha}\right)\right\}c_{i\beta}+\tilde{\mu}_\alpha c_{i\alpha}+\frac{\partial J_{\partial\mathcal{O}}}{\partial f_i}=0.\label{eq:afboundaryu}
\end{equation}
For clarity, we substitute specific values for $i$ in this explanation.
We consider the case where $\Gamma_u$ is on the bottom boundary ($y=0$) using the D2Q9 model, as depicted in Fig.~\ref{fig:boundary_condition}.

To satisfy Eq.~\eqref{eq:afboundaryu}, we need to appropriately assign values to the adjoint variables $\tilde{f}_4$, $\tilde{f}_7$, $\tilde{f}_8$, $\tilde{\mu}_x$ and $\tilde{\mu}_y$, and ensure that the coefficients of $\delta f_2$, $\delta f_4$, $\delta f_5$, $\delta f_6$, $\delta f_7$ and $\delta f_8$ become zero.
However, there are 6 equations against 5 variables and the state is so called over constraint, so they cannot be solved.
Hereupon, we set values to $\tilde{f}_4$, $\tilde{f}_7$, $\tilde{f}_8$, $\tilde{\mu}_x$ and $\tilde{\mu}_y$ with using generalized inverse.
$\bm{M}$, $\bm{N}$, $\tilde{\bm{f}}$, $\bar{\tilde{\bm{f}}}$, $\bm{A}$ and $\bm{B}$ are defined as follows:
\begin{equation}
    \bm{M}=\left(\begin{array}{cccc}
            n_\alpha c_{0\alpha}-2Aw_0n_\alpha c_{0\alpha}c_{0\beta}c_{0\beta} & -2Aw_1n_\alpha c_{1\alpha}c_{1\beta}c_{1\beta}                     & \dots  & -2Aw_8n_\alpha c_{8\alpha}c_{8\beta}c_{8\beta}                      \\
            -2Aw_0n_\alpha c_{0\alpha}c_{0\beta}c_{0\beta}                     & n_\alpha c_{1\alpha}-2Aw_1n_\alpha c_{1\alpha}c_{1\beta}c_{1\beta} & \dots  & -2Aw_8n_\alpha c_{8\alpha}c_{8\beta}c_{8\beta}                      \\
            \vdots                                                             & \vdots                                                             & \ddots & \vdots                                                              \\
            -2Aw_0n_\alpha c_{0\alpha}c_{0\beta}c_{0\beta}                     & -2Aw_1n_\alpha c_{1\alpha}c_{1\beta}c_{1\beta}                     & \dots  & n_\alpha c_{8\alpha}-2Aw_8n_\alpha c_{8\alpha}c_{8\beta}c_{8\beta}
        \end{array}\right),
\end{equation}
\begin{equation}
    \bm{N}=\left(\begin{array}{ccccccccc}
            c_{0x} & c_{1x} & c_{2x} & c_{3x} & c_{4x} & c_{5x} & c_{6x} & c_{7x} & c_{8x} \\
            c_{0y} & c_{1y} & c_{2y} & c_{3y} & c_{4y} & c_{5y} & c_{6y} & c_{7y} & c_{8y}
        \end{array}\right)^T,
\end{equation}
\begin{equation}
    \bm{\tilde{f}}=\left(\begin{array}{ccc}
            \tilde{f}_4 & \tilde{f}_7 & \tilde{f}_8
        \end{array}\right)^T,
    \bm{\bar{\tilde{f}}}=\left(\begin{array}{cccccc}
            \tilde{f}_0 & \tilde{f}_1 & \tilde{f}_2 & \tilde{f}_3 & \tilde{f}_5 & \tilde{f}_6
        \end{array}\right)^T,
\end{equation}
\begin{equation}
    \bm{A}=\left(\begin{array}{ccccccccc}
            0 & 0 & 0 & 0 & 1 & 0 & 0 & 0 & 0  \\
            0 & 0 & 0 & 0 & 0 & 0 & 0 & 1 & 0  \\
            0 & 0 & 0 & 0 & 0 & 0 & 0 & 0 & 1
        \end{array}\right)^T,
    \bm{B}=\left(\begin{array}{ccccccccc}
            1 & 0 & 0 & 0 & 0 & 0 & 0 & 0 & 0  \\
            0 & 1 & 0 & 0 & 0 & 0 & 0 & 0 & 0  \\
            0 & 0 & 1 & 0 & 0 & 0 & 0 & 0 & 0  \\
            0 & 0 & 0 & 1 & 0 & 0 & 0 & 0 & 0  \\
            0 & 0 & 0 & 0 & 0 & 1 & 0 & 0 & 0  \\
            0 & 0 & 0 & 0 & 0 & 0 & 1 & 0 & 0
        \end{array}\right)^T,
\end{equation}
then equations to be solved becomes 
\begin{equation}
    \left(\begin{array}{c|c}
            \bm{M}\bm{A} & \bm{N}
        \end{array}\right)
    \left(\begin{array}{c}
            \bm{\tilde{f}} \\
            \hline
            \bm{\tilde{\mu}}
        \end{array}\right)
    =-\bm{M}\bm{B}\bm{\bar{\tilde{f}}},
\end{equation}
and based on the generalized inverse, they are deformed as 
\begin{equation}
    \left(\begin{array}{c}
            \bm{A}^T\bm{M}^T \\
            \hline
            \bm{N}^T
        \end{array}\right)
    \left(\begin{array}{c|c}
            \bm{M}\bm{A} & \bm{N}
        \end{array}\right)
    \left(\begin{array}{c}
            \bm{\tilde{f}} \\
            \hline
            \bm{\tilde{\mu}}
        \end{array}\right)=-
    \left(\begin{array}{c}
            \bm{A}^T\bm{M}^T \\
            \hline
            \bm{N}^T
        \end{array}\right)
    \bm{M}\bm{B}\bm{\bar{\tilde{f}}}.
\end{equation}
To solve them, we obtain the adjoint variables expressed as
\begin{equation}
    \bm{\tilde{f}}=\left(\begin{array}{c}
            \frac{\tilde{f}_2+\tilde{f}_5+\tilde{f}_6}{3}    \\
            \frac{4\tilde{f}_2+7\tilde{f}_5+\tilde{f}_6}{12} \\
            \frac{4\tilde{f}_2+\tilde{f}_5+7\tilde{f}_6}{12}
        \end{array}\right),
    \bm{\tilde{\mu}}=\left(\begin{array}{c}
            -\frac{\left(A-3\right)\left(\tilde{f}_5-\tilde{f}_6\right)}{12} \\
            -\frac{2\left(A-1\right)\tilde{f}_2+\left(A-2\right)\left(\tilde{f}_5+\tilde{f}_6\right)}{6}
        \end{array}\right).
    \label{eq:ftilde}
\end{equation}
When solving the adjoint equations, $\tilde{\rho}$, $\tilde{u}_\alpha$ and $\tilde{s}_{\alpha\beta}$ from Eq.~\eqref{eq:ftilde} are applied directly at the boundary.

Considering meanings of $\bm{\tilde{f}}$ derived in the above, for the over constraint equations $\bm{A}\bm{x}=\bm{b}$ the solutions $\bm{x}$ given by generalized inverse make $\frac{1}{2}\left\|\bm{A}\bm{x}-\bm{b}\right\|_2^2$ minimize. 
Thus, this method yields $\bm{\tilde{f}}$, which minimizes the norm of $\bm{r}$, defined as:
\begin{equation}
    r_i:=n_\alpha\left\{\tilde{f}_i\delta_{\alpha\beta}-A\left(\tilde{s}_{\alpha\beta}+\tilde{s}_{\beta\alpha}\right)\right\}c_{i\beta}+\tilde{\mu}_\alpha c_{i\alpha}, \label{eq:r_definition}
\end{equation}
where $r_i$ represents the coefficient for $\delta f_i$ on the inlet boundary $\Gamma_u$.

\subsection{Sensitivity verification}
\label{sec26}

In this section, we verify sensitivities of non-thermal/thermal fluid problems based on ALKS by comparing with sensitivities based on finite difference approximation. 

\begin{figure}[t]
    \centering
    \begin{minipage}[b]{0.33\columnwidth}
        \centering
        \includegraphics[width=0.9\columnwidth]{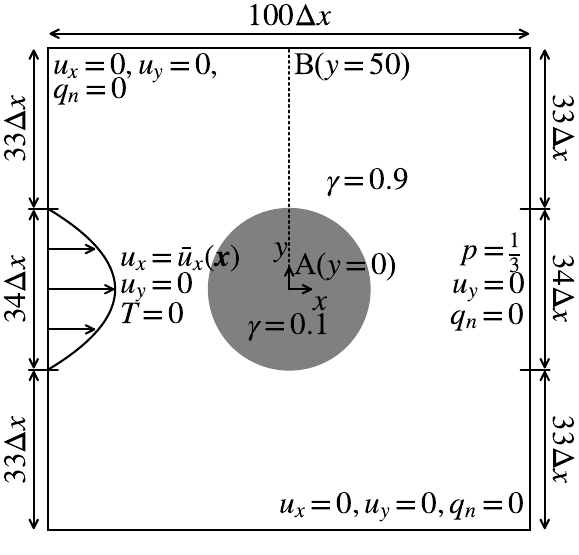}
        \subcaption{Design setting}
        \label{fig:scsens}
    \end{minipage}
    \begin{minipage}[b]{0.33\columnwidth}
        \centering
        \includegraphics[width=0.9\columnwidth]{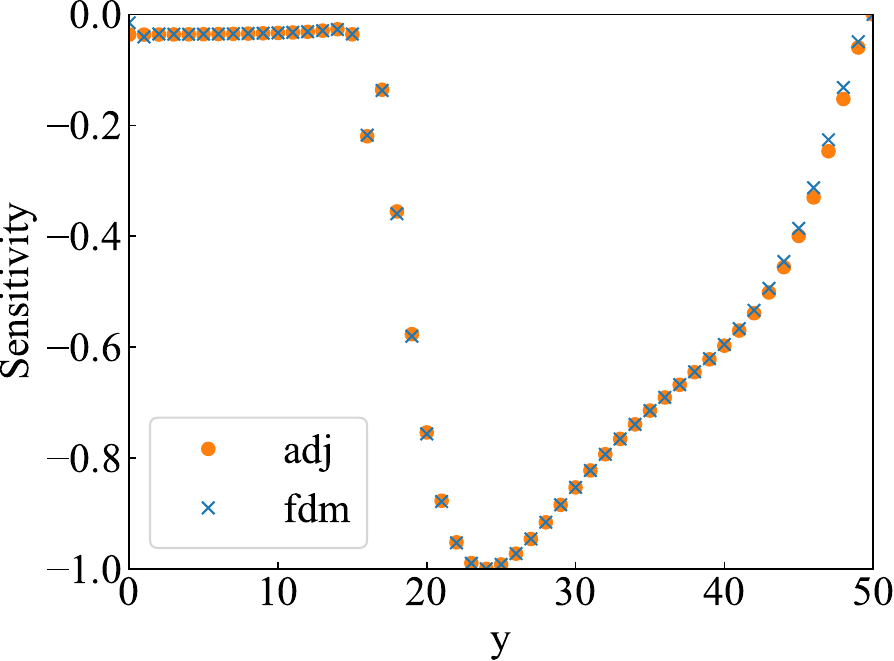}
        \subcaption{Result of non-thermal fluid problem}
        \label{fig:nssens}
    \end{minipage}
    \begin{minipage}[b]{0.33\columnwidth}
        \centering
        \includegraphics[width=0.9\columnwidth]{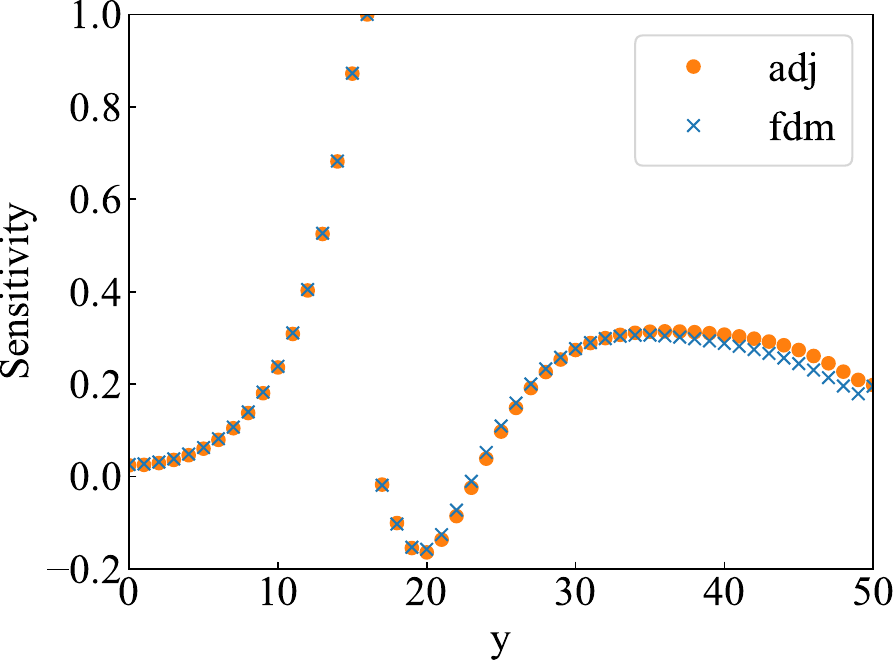}
        \subcaption{Result of forced convection problem}
        \label{fig:nsadsens}
    \end{minipage}
    \caption{Design setting and results of sensitivity verification for non-thermal fluid problem and forced convection problem}
\end{figure}

First, as non-thermal fluid problem, we consider the pressure drop in a square region which is surrounded by walls depicted in Fig.~\ref{fig:scsens}.
Fluid flows into the region from part of the left boundary with a parabolic profile expressed as:
\begin{equation}
    \bar{u}_x\left(\bm{x}\right)=u_0\left(y-33\Delta x\right)\left(y-67\Delta x\right), \label{eq:inlet_flow}
\end{equation}
and exits from a part of the right boundary.
Analysis domain consists of $100\Delta x\times 100\Delta x$ grids and design valiable $\gamma$ is set to $0.1$ in the center circle whose radius is $100/6\Delta x$ and in the other region $\gamma$ is set to $0.9$.
The kinematic viscosity $\nu$ is set to $0.1$, and the maximum fluid velocity at the inlet boundary $u_0$ is adjusted to satisfy a Reynolds number of $Re=1$ for a characteristic length of $100/3\Delta x$.
The finite difference approximation for sensitivity is calculated as follows:    
\begin{equation}
    J_{\text{FD}n}=\frac{J_n\left(\gamma+\epsilon\right)-J_n\left(\gamma-\epsilon\right)}{2\epsilon},\label{eq:sensfdm}
\end{equation}
where $\epsilon$ is set to $10^{-3}$.
In the non-thermal fluid problem, the boundary conditions for thermal variables such as $T$ and $q_n$ are ignored.
The sensitivities are compared along the dashed line in the Fig.~\ref{fig:scsens} and the result is in the Fig.~\ref{fig:nssens}.
``adj'' and ``fdm'' in the Fig.~\ref{fig:nssens} mean the sensitivities based on ALKS and FDM respectively and there is strong agreement between them.

Subsequently, the sensitivity of the forced convection problem is verified using the same method as for the non-thermal fluid problem.
The dimensions and boundary conditions are as in the Fig.~\ref{fig:nsadsens}, and temperature is prescribed on the inlet boundary and on the other boundaries are adiabatic.
The boundary conditions for the pressure and the flow velocity are the same with the non-thermal fluid problem. 
The result is in the Fig.~\ref{fig:nsadsens}, and the sensitivity based on ALKS is fairly consistent with that of FDM, although some errors are observed near the boundary. 
In addition to the effects of approximations introduced in the boundary conditions of the adjoint equations (as discussed in Section~\ref{sec25}), according to previous studies~\cite{Nadarajah2000comparison}, the continuous adjoint approach is generally less accurate than the discrete adjoint approach.

\begin{figure}[t]
    \centering
    \begin{minipage}[b]{0.49\columnwidth}
        \centering
        \includegraphics[width=0.9\columnwidth]{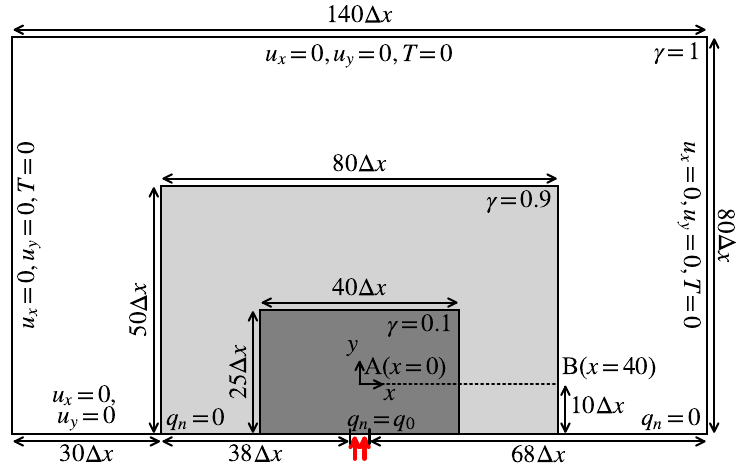}
        \subcaption{Design setting}
        \label{fig:scnsadncsens}
    \end{minipage}
    \begin{minipage}[b]{0.49\columnwidth}
        \centering
        \includegraphics[width=0.6\columnwidth]{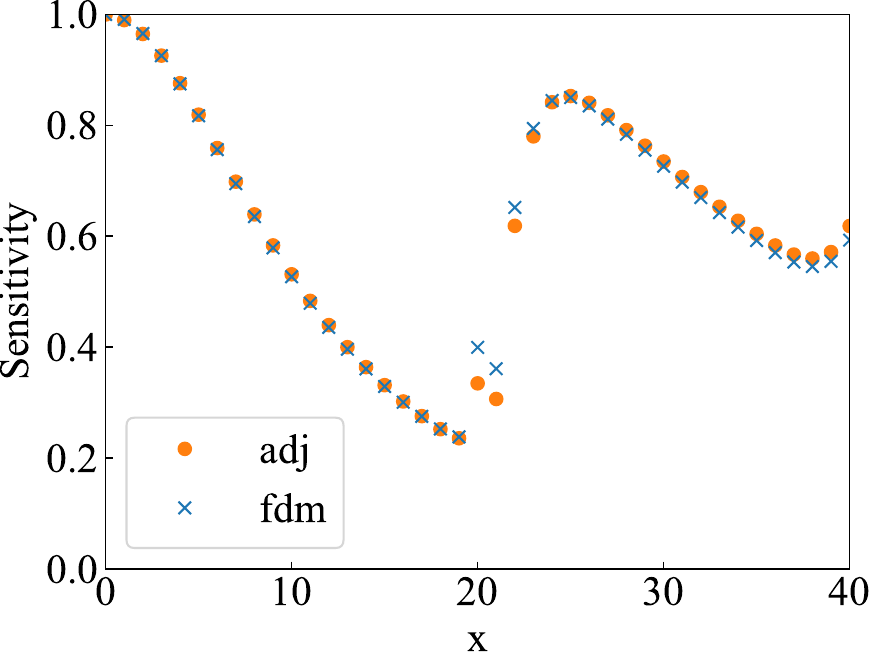}
        \subcaption{Result of natural convection problem}
        \label{fig:nsadncsens}
    \end{minipage}
    \caption{Design setting and result of sensitivity verification for natural convection problem}
\end{figure}

Finally, we also verify sensitivities of natural convection problems.
The analysis domain consists of $140\Delta x\times80\Delta x$ grids and is surrounded by walls as shown in Fig.~\ref{fig:scnsadncsens}.
On the right, left and top boundary temperature prescribed boundary is applied, on the other hand, on the bottom boundary heating and adiabatic boundaries are applied at center $4\Delta x$ part and the other, respectively.
The design domain consists of the central $80\Delta x\times50\Delta x$ grid, and in its central $40\Delta x\times25\Delta x$ part $\gamma$ is set to $0.1$ while in the other part $\gamma$ is set to $0.9$, furthermore the non-design domain is fluid, in other words $\gamma$ is set to $1.0$.
The sensitivities are compared along the dashed line in the Fig.~\ref{fig:scnsadncsens} and the result is in the Fig.~\ref{fig:nsadncsens}.
Although some errors are observed around the interface between the regions which $\gamma=0.1$ and $\gamma=0.9$, the sensitivity based on ALKS is fairly consistent with that of FDM.

\subsection{Optimization procedure}
\label{sec27}

The numerical procedure of the proposed method is as follows:

\begin{enumerate}
    \item[{\it Step 1}] The design variable field $\gamma$ is initialized on the analysis domain discretized using square lattice mesh.
    \item[{\it Step 2}] The state variable fields $\rho$, $u_\alpha$, $T$ are computed using the LKS.
    \item[{\it Step 3}] The values of the objective functional $J$ and the constraint functional $G$ are computed.
    \item[{\it Step 4}] The adjoint variable fields $\tilde{\rho}$, $\tilde{u}_\alpha$, $\tilde{s}_{\alpha\beta}$, $\tilde{T}$, $\tilde{q}_\alpha$ are computed using the ALKS.
    \item[{\it Step 5}] The design sensitivities $J^\prime$ and $G^\prime$ are computed from the state and adjoint variable fields.
    \item[{\it Step 6}] The design variable field $\gamma$ is updated using the method of moving asymptotes (MMA)~\cite{svanberg1987method}.
    \item[{\it Step 7}] The procedure returns to the {\it Step 2} of the iteration loop until the convergence criterion is met.
\end{enumerate}

\section{Numerical examples}
\label{sec3}

In this section, we present numerical examples covering both thermal and non-thermal cases, as well as steady and unsteady fluid flows.
The validity of the optimized configurations is evaluated through physical interpretation and by referencing previous studies.
Additionally, we examine and compare the memory requirements of each example using the proposed LKS method against those of the LBM approach.

\subsection{Pipe bend}
\label{sec31}

\begin{figure}[t]
    \centering
    \begin{minipage}[b]{0.33\columnwidth}
        \centering
        \includegraphics[width=0.9\columnwidth]{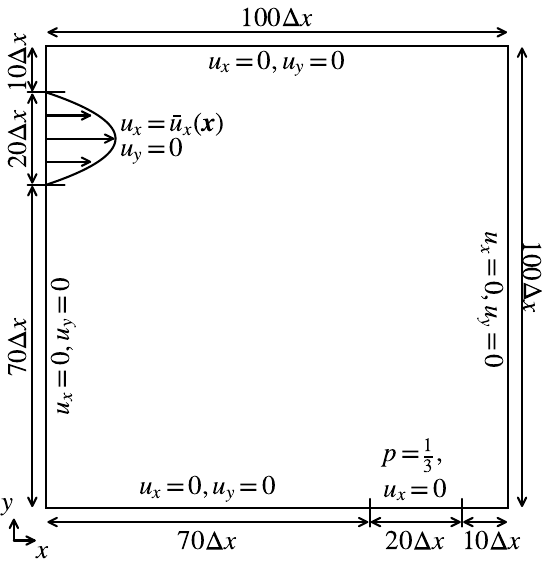}
        \subcaption{Design setting}
        \label{fig:design_setting_pipebend}
    \end{minipage}
    \begin{minipage}[b]{0.33\columnwidth}
        \centering
        \includegraphics[width=0.9\columnwidth]{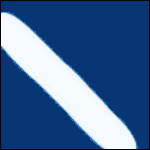}
        \subcaption{optimized configuration for $Re=1$}
        \label{fig:optimized_shape_pipebend_1}
    \end{minipage}
    \begin{minipage}[b]{0.33\columnwidth}
        \centering
        \includegraphics[width=0.9\columnwidth]{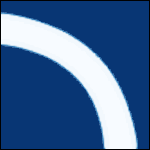}
        \subcaption{optimized configuration for $Re=100$}
        \label{fig:optimized_shape_pipebend_2}
    \end{minipage}
    \caption{Design setting and optimized configurations of pipe bend problem}
\end{figure}

We first address the design of a bent flow path as a numerical example of a steady-state non-thermal fluid problem.
As shown in Fig.~\ref{fig:design_setting_pipebend}, fluid flows into the analysis domain from the top left with a parabolic profile expressed as:
\begin{equation}
    \bar{u}_x\left(\bm{x}\right)=u_0\left(y-70\Delta x\right)\left(y-90\Delta x\right),
\end{equation}
and exits from the bottom right.
The objective functional of this problem is to minimize the static pressure drop expressed as:
\begin{equation}
    J=\int_{\partial\mathcal{O}_\text{in}}pd\Gamma-\int_{\partial\mathcal{O}_\text{out}}pd\Gamma. \label{eq:objective_pipebend}
\end{equation}
Here, $\partial\mathcal{O}_\text{in}$ and $\partial\mathcal{O}_\text{out}$ represent the inlet and outlet boundaries, respectively.
Additionally, in this problem, the flow path area constraint expressed as: 
\begin{equation}
    G=\int_\mathcal{O}\gamma d\Omega\leq V_\text{max}\int_\mathcal{O}d\Omega \label{eq:constraint_pipebend}
\end{equation}
is applied with its maximum value $V_\text{max}$ set to 0.25, corresponding to 25\% of the analysis domain.

Fig.~\ref{fig:optimized_shape_pipebend_1} and \ref{fig:optimized_shape_pipebend_2} present the optimized configurations for $Re=1$ and $Re=100$, respectively.
The Reynolds number for this problem is defined as:
\begin{equation}
    Re=\frac{Lu_0}{\nu}, \label{eq:reynolds_pipebend}
\end{equation}
where the characteristic length $L$ set to $20\Delta x$, corresponding to the inlet width, and the kinematic viscosity $\nu$ is set to $0.1$.
For low Reynolds numbers, the shape of the flow path resembles a straight line connecting the inlet and the outlet, whereas for high Reynolds numbers, the shape forms an arc.
This effect can be attributed to the stronger influence of flow inertia at higher Reynolds numbers.
Moreover, a similar trend regarding the influence of the Reynolds number has been observed in previous research~\cite{gersborg2005topology}.

The estimated memory required by ALKS for this example is shown in Table~\ref{tab:pipebend_memory_usage_ALKS}, while Table~\ref{tab:pipebend_memory_usage_ALBM} presents the memory required by ALBM under identical conditions.
Consequently, ALKS reduces memory usage by approximately 42\% compared to ALBM for sensitivity analysis.
\begin{table}[t]
    \centering
    \begin{minipage}{\textwidth}
        \centering
        \caption{Memory requirements for pipe bend problem using ALKS}
        \begin{tabular}{|p{0.45\linewidth}|c|c|} \hline
            Variable name                                                                                                                                                                                                                                                  & Array length            & Total amount                                    \\ \hline
            $\partial u_x/\partial x$,$\partial u_y/\partial x$,$\partial u_x/\partial y$,$\partial u_y/\partial y$,\newline$\partial\tilde{s}_{xx}/\partial x$,$\partial\tilde{s}_{xy}\partial x$,$\partial\tilde{s}_{yx}/\partial y$,$\partial\tilde{s}_{yy}/\partial y$ & $n_x\times n_y$         & $101\times 101\times 8\times 8=653$kB           \\ \hline
            $\rho$,$u_x$,$u_y$,$\tilde{\rho}$,$\tilde{u}_x$,$\tilde{u}_y$,$\tilde{s}_{xx}$,$\tilde{s}_{xy}$,$\tilde{s}_{yx}$,$\tilde{s}_{yy}$                                                                                                                              & $n_x\times n_y\times 2$ & $101\times 101\times 2\times 10\times 8=1.63$MB \\ \hline
        \end{tabular}
        \label{tab:pipebend_memory_usage_ALKS}
    \end{minipage}
    \begin{minipage}{\textwidth}
        \centering
        \caption{Memory requirements for pipe bend problem using ALBM}
        \begin{tabular}{|c|c|c|} \hline
            Variable name                                                                         & Array length            & Total amount                                    \\ \hline
            $\rho$,$\tilde{m}_x$,$\tilde{m}_y$,$\tilde{\rho}$                                     & $n_x\times n_y$         & $101\times 101\times 4\times 8=326$kB           \\ \hline
            $u_x$,$u_y$,$f_0\dots f_8$,$\tilde{u}_x$,$\tilde{u}_y$,$\tilde{f}_0\dots \tilde{f}_8$ & $n_x\times n_y\times 2$ & $101\times 101\times 2\times 22\times 8=3.59$MB \\ \hline
        \end{tabular}
        \label{tab:pipebend_memory_usage_ALBM}
    \end{minipage}
\end{table}

\subsection{Periodic double pipe}
\label{sec32}

\begin{figure}[t]
    \centering
    \begin{minipage}[b]{0.49\columnwidth}
        \centering
        \includegraphics[width=0.65\columnwidth]{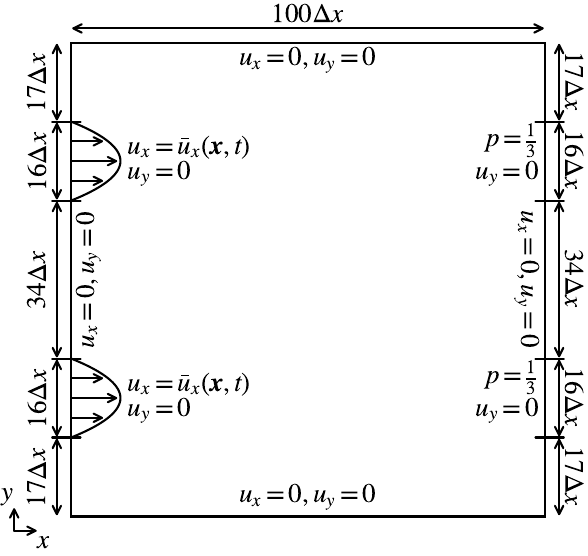}
        \subcaption{Design setting}
        \label{fig:design_setting_doublepipe}
    \end{minipage}
    \begin{minipage}[b]{0.49\columnwidth}
        \centering
        \includegraphics[width=0.6\columnwidth]{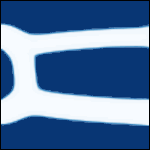}
        \subcaption{optimized configuration}
        \label{fig:optimized_shape_doublepipe}
    \end{minipage}
    \caption{Design setting and optimized configuration of periodic double pipe problem}
\end{figure}

Next, as the unsteady non-thermal fluid problem, we address the double pipe flow path design.
As shown in Fig.~\ref{fig:design_setting_doublepipe}, fluid flows into the analysis domain from the left top and bottom with parabolic profiles such as:
\begin{align}
    \bar{u}_x\left(\bm{x},t\right)=\left\{
    \begin{array}{cl}
        u_0\cos\left(\frac{2\pi t}{N_t}\right)\left(y-17\Delta x\right)\left(y-33\Delta x\right) & \mbox{if $17\Delta x\leq y\leq 33\Delta x$}  \\
        u_0\sin\left(\frac{2\pi t}{N_t}\right)\left(y-57\Delta x\right)\left(y-83\Delta x\right) & \mbox{if $57\Delta x\leq y\leq 83\Delta x$}, \\
    \end{array}
    \right.
\end{align}
and exits from the right top and bottom boundary.
The objective functional of this problem is described as minimizing the loss of the energy expressed as: 
\begin{equation}
    J=\int_\mathcal{I}\int_{\partial\mathcal{O}}-n_\alpha u_\alpha\left(p+\frac{1}{2}u_\beta^2\right)d\Gamma dt. \label{eq:objective_doublepipe}
\end{equation}
Similar to the bended flow path design problem in Section \ref{sec31}, a flow path area constraint is applied in this problem, with the maximum value set to 0.33, corresponding to 33\% of the analysis domain.
The Reynolds number $Re$ is set to $1$, and the characteristic length $L$ is set to $100\Delta x/6$ which is equal to the width of the inlet and the kinematic viscosity $\nu$ is set to $0.1$.

\begin{figure}[t]
    \centering
    \begin{minipage}[b]{0.49\columnwidth}
        \centering
        \includegraphics[width=0.6\columnwidth]{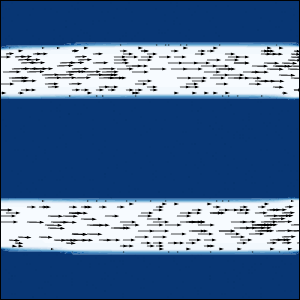}
        \subcaption{``push-push''}
        \label{fig:optimized_shape_doublepipe_push_push}
    \end{minipage}
    \begin{minipage}[b]{0.49\columnwidth}
        \centering
        \includegraphics[width=0.6\columnwidth]{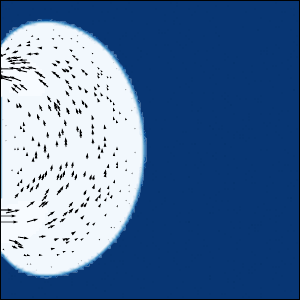}
        \subcaption{``push-pull''}
        \label{fig:optimized_shape_doublepipe_push_pull}
    \end{minipage}
    \caption{Optimized configurations for ``push-push'' and ``push-pull'' cases}
\end{figure}

\begin{figure}[!ht]
    \centering
    \begin{minipage}[b]{0.33\columnwidth}
        \centering
        \includegraphics[width=0.9\columnwidth]{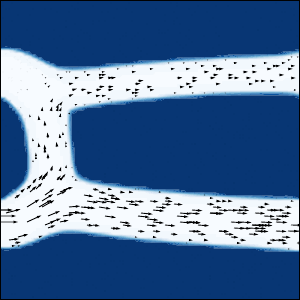}
        \subcaption{$t=0$}
    \end{minipage}
    \begin{minipage}[b]{0.33\columnwidth}
        \centering
        \includegraphics[width=0.9\columnwidth]{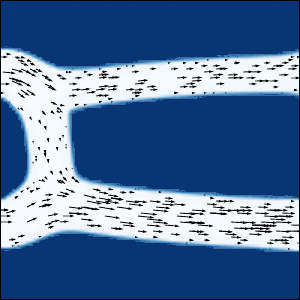}
        \subcaption{$t=0.125$}
    \end{minipage}
    \begin{minipage}[b]{0.33\columnwidth}
        \centering
        \includegraphics[width=0.9\columnwidth]{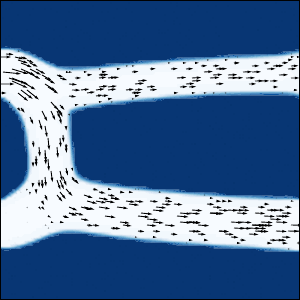}
        \subcaption{$t=0.25$}
    \end{minipage}
    \\
    \begin{minipage}[b]{0.33\columnwidth}
        \centering
        \includegraphics[width=0.9\columnwidth]{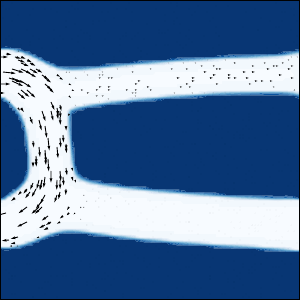}
        \subcaption{$t=0.375$}
    \end{minipage}
    \begin{minipage}[b]{0.33\columnwidth}
        \centering
        \includegraphics[width=0.9\columnwidth]{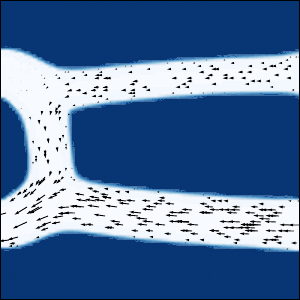}
        \subcaption{$t=0.5$}
    \end{minipage}
    \begin{minipage}[b]{0.33\columnwidth}
        \centering
        \includegraphics[width=0.9\columnwidth]{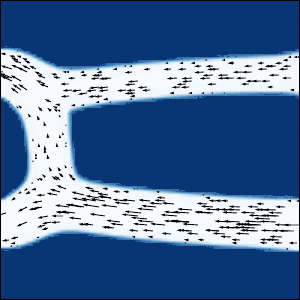}
        \subcaption{$t=0.625$}
    \end{minipage}
    \\
    \begin{minipage}[b]{0.33\columnwidth}
        \centering
        \includegraphics[width=0.9\columnwidth]{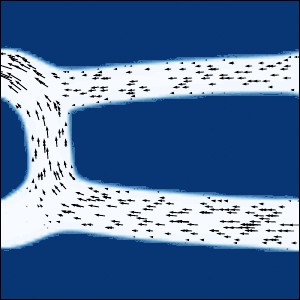}
        \subcaption{$t=0.75$}
    \end{minipage}
    \begin{minipage}[b]{0.33\columnwidth}
        \centering
        \includegraphics[width=0.9\columnwidth]{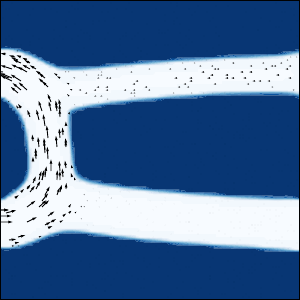}
        \subcaption{$t=0.875$}
    \end{minipage}
    \begin{minipage}[b]{0.33\columnwidth}
        \centering
        \includegraphics[width=0.9\columnwidth]{assets_example2_doublepipe_case17_velocity_0.pdf}
        \subcaption{$t=1.0$}
    \end{minipage}
    \caption{Transition of the velocity distribution passing through the optimized configuration}
    \label{fig:velocity_doublepipe}
\end{figure}

Fig.~\ref{fig:optimized_shape_doublepipe} shows the optimized configuration for the double pipe flow path design.
The design takes an ``H'' shape, resembling that presented in the previous study~\cite{DENG20116688}.
Although researchers have previously explained the reason for the optimized ``H'' shape, we reiterate it here for clarity.
At the inlets on the top left and bottom left, the flow alternates between ``push-pull'', ``push-push'', ``pull-push'' and ``pull-pull'', resulting in an optimized configuration that is a compromise between ``push-push'' and ``push-pull''.
The optimized configurations for each case are depicted in Figs.~\ref{fig:optimized_shape_doublepipe_push_push} and \ref{fig:optimized_shape_doublepipe_push_pull}, furthermore, the flow in the optimized "H" shape changes sequentially as shown in Fig.~\ref{fig:velocity_doublepipe} (the figures at $t=1.0$ and $t=0$ are the same images).

Table~\ref{tab:doublepipe_memory_usage_ALKS} provides the estimated memory usage for this example using the proposed method, while the memory required by ALBM under the same conditions is estimated as Table~\ref{tab:doublepipe_memory_usage_ALBM}.
Notably, the memory required for sensitivity analysis using ALKS is approximately equal to that of ALBM. 
This is primarily because, in unsteady non-thermal fluid simulations, the majority of memory is allocated to storing the velocity field across each time step, a requirement shared by both ALKS and ALBM.
\begin{table}[t]
    \centering
    \begin{minipage}{\textwidth}
        \centering
        \caption{Memory requirements for periodic double pipe problem using ALKS}
        \begin{tabular}{|p{0.45\linewidth}|c|c|} \hline
            Variable name                                                                                                                                                                                                                                                  & Array length              & Total amount                                       \\ \hline
            $\partial u_x/\partial x$,$\partial u_y/\partial x$,$\partial u_x/\partial y$,$\partial u_y/\partial y$,\newline$\partial\tilde{s}_{xx}/\partial x$,$\partial\tilde{s}_{xy}\partial x$,$\partial\tilde{s}_{yx}/\partial y$,$\partial\tilde{s}_{yy}/\partial y$ & $n_x\times n_y$           & $101\times 101\times 8\times 8=653$kB              \\ \hline
            $\rho$,$\tilde{\rho}$,$\tilde{u}_x$,$\tilde{u}_y$,$\tilde{s}_{xx}$,$\tilde{s}_{xy}$,$\tilde{s}_{yx}$,$\tilde{s}_{yy}$                                                                                                                                          & $n_x\times n_y\times 2$   & $101\times 101\times 2\times 8\times 8=1.31$MB     \\ \hline
            $u_x$,$u_y$                                                                                                                                                                                                                                                    & $n_x\times n_y\times n_t$ & $101\times 101\times 20000\times 2\times 8=3.26$GB \\\hline
        \end{tabular}
        \label{tab:doublepipe_memory_usage_ALKS}
    \end{minipage}
    \begin{minipage}{\textwidth}
        \centering
        \caption{Memory requirements for periodic double pipe problem using ALBM}
        \begin{tabular}{|c|c|c|} \hline
            Variable name                                                             & Array length              & Total amount                                       \\ \hline
            $\rho$,$\tilde{m}_x$,$\tilde{m}_y$,$\tilde{\rho}$                         & $n_x\times n_y$           & $101\times 101\times 4\times 8=326$kB              \\ \hline
            $f_0\dots f_8$,$\tilde{u}_x$,$\tilde{u}_y$,$\tilde{f}_0\dots \tilde{f}_8$ & $n_x\times n_y\times 2$   & $101\times 101\times 2\times 20\times 8=3.26$MB    \\ \hline
            $u_x$,$u_y$                                                               & $n_x\times n_y\times n_t$ & $101\times 101\times 20000\times 2\times 8=3.26$GB \\\hline
        \end{tabular}
        \label{tab:doublepipe_memory_usage_ALBM}
    \end{minipage}
\end{table}

\subsection{Heat exchanger}
\label{sec33}

\begin{figure}[t]
    \centering
    \begin{minipage}[b]{0.49\columnwidth}
        \centering
        \includegraphics[width=0.6\columnwidth]{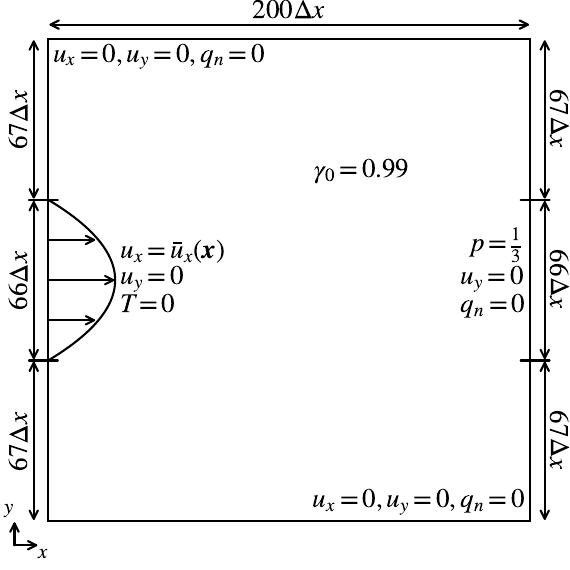}
        \subcaption{Design setting}
        \label{fig:design_setting_heatexchanger}
    \end{minipage}
    \begin{minipage}[b]{0.49\columnwidth}
        \centering
        \includegraphics[width=0.6\columnwidth]{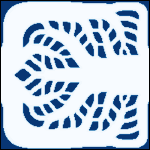}
        \subcaption{optimized configuration}
        \label{fig:optimized_shape_heatexchanger}
    \end{minipage}
    \\
    \begin{minipage}[b]{0.49\columnwidth}
        \centering
        \includegraphics[width=0.65\columnwidth]{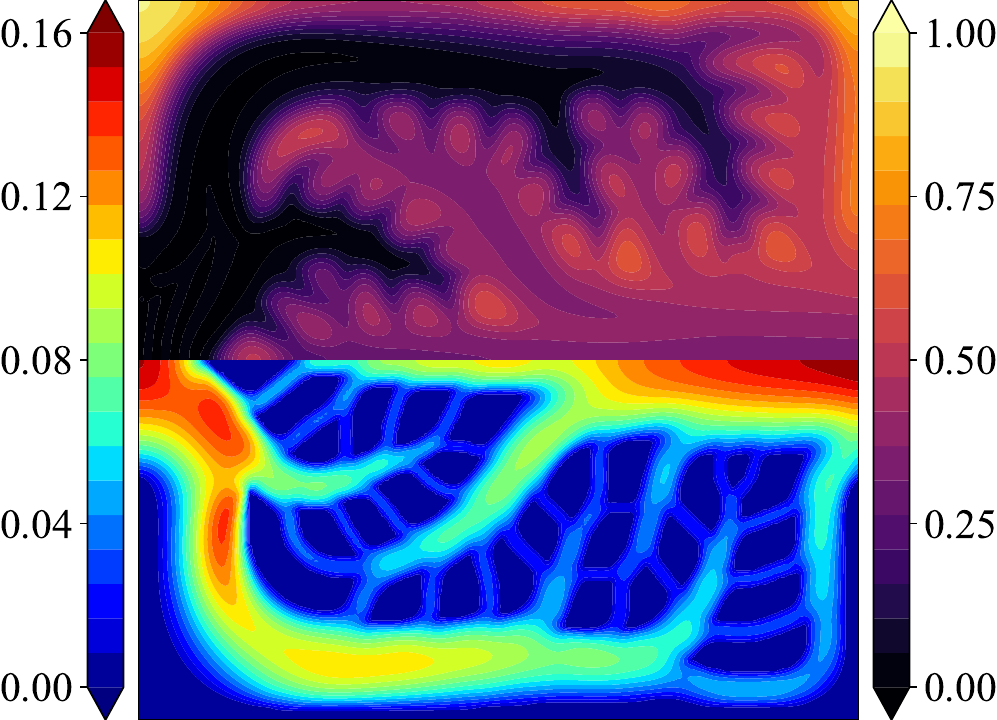}
        \subcaption{Velocity and temperature distributions}
        \label{fig:velocity_temperature_heatexchanger}
    \end{minipage}
    \begin{minipage}[b]{0.49\columnwidth}
        \centering
        \includegraphics[width=0.7\columnwidth]{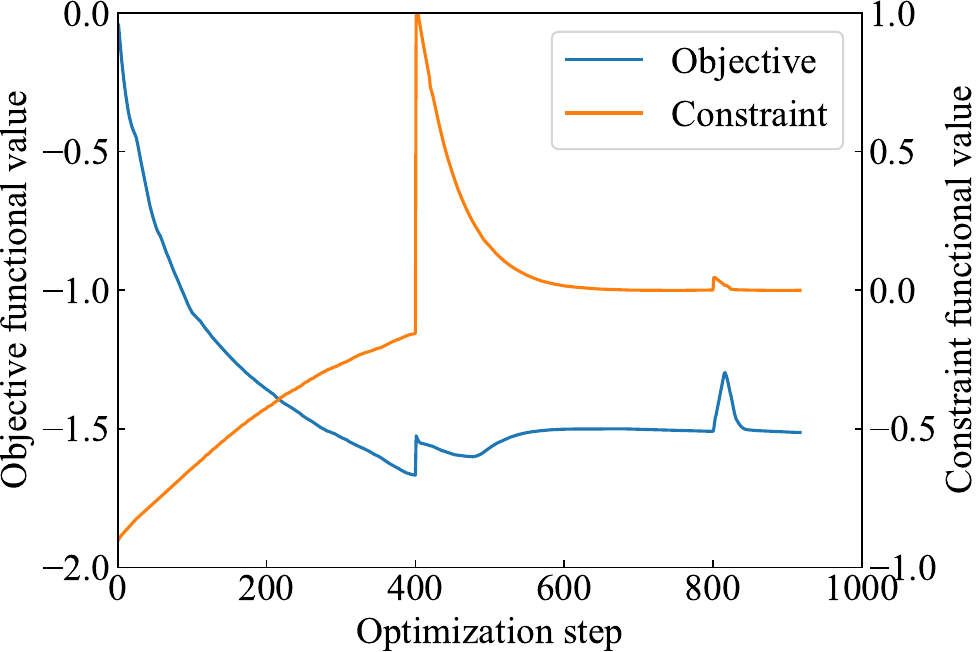}
        \subcaption{History plot of objective and constraint functional values}
        \label{fig:historical_values_heatexchanger}
    \end{minipage}
    \caption{Design setting, optimized configuration, velocity and temperature distribution, and history plot of objective and constraint functional values of heat exchanger problem}
\end{figure}

Next, for the steady-state forced convection problem, we address the design of a heat exchanger.
As shown in Fig.~\ref{fig:design_setting_heatexchanger}, the coolant enters into the analysis domain from the left center with a parabolic profile expressed as:
\begin{equation}
    \bar{u}_x\left(\bm{x}\right)=u_0\left(y-67\Delta x\right)\left(y-133\Delta x\right),
\end{equation}
and exits from the right center boundary.
The objective functional of this problem is described as maximizing the heat exchange amount expressed as:
\begin{equation}
    J=-\int_\mathcal{O}\beta_\gamma\left(1-T\right)d\Omega. \label{eq:objective_heatexchanger}
\end{equation}
In this problem, the source term of the energy equation (Eq.~\eqref{eq:energy}) is as follows:
\begin{equation}
    Q=\beta_\gamma\left(1-T\right).\label{eq:heatsource}
\end{equation}
This refers to a fictitious source generating heat only in the solid region and $\beta_\gamma$ in Eq.~\eqref{eq:heatsource} is defined as:   
\begin{equation}
    \beta_\gamma=\bar{\beta}\frac{q_\beta\left(1-\gamma\right)}{q_\beta+\gamma}.\label{eq:betagamma}
\end{equation}
For the heat exchanger design problem the maximum value of the static pressure drop is constrained and it is expressed as follows:
\begin{equation}
    G=\int_{\partial\mathcal{O}_\text{in}}pd\Gamma-\int_{\partial\mathcal{O}_\text{out}}pd\Gamma\leq\eta\Delta p_0, \label{eq:constraint_heatexchanger}
\end{equation}
where $\Delta p_0$ is the static pressure drop at first step of the optimization loop.
The Reynolds number $Re$ is set to $100$, with the characteristic length $L$ equal to $200\Delta x/3$, matching the inlet width, and the kinematic viscosity $\nu$ is set to $0.1$.
Besides, the Prandtl number $Pr$ is set to $6$.

The optimized configuration is shown in Fig.~\ref{fig:optimized_shape_heatexchanger} and the upper half of Fig.~\ref{fig:velocity_temperature_heatexchanger} shows the temperature distribution against the lower half part shows the flow velocity distribution.
The shape of the flow path is designed to maximize heat exchange by increasing the interface between the solid and fluid regions, and this trend is consistent with previous studies~\cite{Matsumori2013}.
The convergence histories of the objective and the constraint functionals values over the optimization iteration are plotted in Fig.~\ref{fig:historical_values_heatexchanger}.
The peak in both functional values in Fig.~\ref{fig:historical_values_heatexchanger} is likely due to the Continuation scheme.
In this case, the parameter value $q_\alpha$ in the Brinkman model is updated every $400$ steps of the optimization iteration and its curve gradually becomes linear.
This is likely the reason why, after updating, the fluid flow is restricted in the grayscale region, causing the heat exchange amount and the static pressure drop increase dramatically.
\begin{figure}[t]
    \centering
    \begin{minipage}[b]{0.33\columnwidth}
        \centering
        \includegraphics[width=0.9\columnwidth]{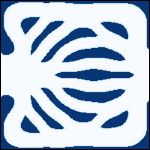}
        \subcaption{$Re=10$, $b_\text{max}=0.1$, $\eta_\text{max}=10$}
        \label{fig:optimized_shape_heatexchanger_2}
    \end{minipage}
    \begin{minipage}[b]{0.33\columnwidth}
        \centering
        \includegraphics[width=0.9\columnwidth]{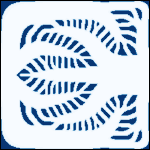}
        \subcaption{$Re=100$, $b_\text{max}=0.3$, $\eta_\text{max}=10$}
        \label{fig:optimized_shape_heatexchanger_3}
    \end{minipage}
    \begin{minipage}[b]{0.33\columnwidth}
        \centering
        \includegraphics[width=0.9\columnwidth]{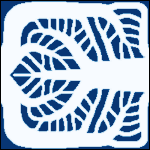}
        \subcaption{$Re=100$, $b_\text{max}=0.1$, $\eta_\text{max}=20$}
        \label{fig:optimized_shape_heatexchanger_4}
    \end{minipage}
    \caption{Optimized configuration for each parameter}
\end{figure}

We also investigate the effects of each parameter for the shape.
First, about the effect of the Reynold number, Fig.~\ref{fig:optimized_shape_heatexchanger_2} shows the optimized configuration for $Re=10$.
When the Reynolds number is low, the flow velocity is deliberately set to a lower value, resulting in a reduced inlet pressure required to drive the fluid during the initial phase of the optimization iterations. 
Consequently, the constraint on the maximum pressure drop becomes more stringent compared to the high Reynolds number scenario, leading to the formation of a simpler channel.
Next, about the effect of the heat generation amount in the solid region, Fig.~\ref{fig:optimized_shape_heatexchanger_3} shows the optimized configuration for $b_\text{max}=0.3$.
As the amount of heat generation increases, increasing the interfaces between the solid and fluid becomes more effective, resulting in finer flow channels.
Finally, about the maximum value of the static pressure drop allowed, Fig.~\ref{fig:optimized_shape_heatexchanger_4} shows the optimized configuration for $\eta_\text{max}=20$.
As the maximum value increases, the flow channel becomes more complex, promoting more heat exchange due to the increased interfaces between the solid and fluid.
Furthermore, the qualitative dependencies of each parameter are consistent with the results in previous research~\cite{Matsumori2013}\cite{YAJI2016355}.

The memory required by ALKS for this example is shown in Table~\ref{tab:heatexchanger_memory_usage_ALKS}, while the memory required by ALBM under the same conditions is shown in Table~\ref{tab:heatexchanger_memory_usage_ALBM}.
Consequently, for the sensitivity analysis, ALKS requires about 53\% less memory than ALBM.
\begin{table}[t]
    \centering
    \begin{minipage}{\textwidth}
        \centering
        \caption{Memory requirements for heat exchanger problem using ALKS}
        \begin{tabular}{|p{0.45\linewidth}|c|c|} \hline
            Variable name                                                                                                                                                                                                                                                                                                                                                                                & Array length            & Total amount                                    \\ \hline
            $\partial u_x/\partial x$,$\partial u_y/\partial x$,$\partial u_x/\partial y$,$\partial u_y/\partial y$, \newline $\partial\tilde{s}_{xx}/\partial x$,$\partial\tilde{s}_{xy}\partial x$,$\partial\tilde{s}_{yx}/\partial y$,$\partial\tilde{s}_{yy}/\partial y$, \newline $\partial T/\partial x$,$\partial T/\partial y$,$\partial\tilde{q}_x/\partial x$,$\partial\tilde{q}_y/\partial y$ & $n_x\times n_y$         & $201\times 201\times 12\times 8=3.88$MB         \\ \hline
            $\rho$,$u_x$,$u_y$,$T$,$q_x$,$q_y$,$\tilde{\rho}$,$\tilde{u}_x$,$\tilde{u}_y$,$\tilde{s}_{xx}$,$\tilde{s}_{xy}$,$\tilde{s}_{yx}$,$\tilde{s}_{yy}$,$\tilde{T}$,$\tilde{q}_x$,$\tilde{q}_y$                                                                                                                                                                                                    & $n_x\times n_y\times 2$ & $201\times 201\times 2\times 16\times 8=10.3$MB \\ \hline
        \end{tabular}
        \label{tab:heatexchanger_memory_usage_ALKS}
    \end{minipage}
    \begin{minipage}{\textwidth}
        \centering
        \caption{Memory requirements for heat exchanger problem using ALBM}
        \begin{tabular}{|c|c|c|} \hline
            Variable name                                                                                                                                                               & Array length            & Total amount                                    \\ \hline
            $\rho$,$T$,$\tilde{m}_x$,$\tilde{m}_y$,$\tilde{\rho}$,$\tilde{T}$                                                                                                           & $n_x\times n_y$         & $201\times 201\times 6\times 8=1.94$MB          \\ \hline
            $u_x$,$u_y$,$f_0\dots f_8$,$\tilde{u}_x$,$\tilde{u}_y$,$\tilde{f}_0\dots\tilde{f}_8$,$q_x$,$q_y$,$g_0\dots g_8$,$\tilde{q}_{x}$,$\tilde{q}_y$,$\tilde{g}_0\dots\tilde{g}_8$ & $n_x\times n_y\times 2$ & $201\times 201\times 2\times 44\times 8=$28.4MB \\ \hline
        \end{tabular}
        \label{tab:heatexchanger_memory_usage_ALBM}
    \end{minipage}
\end{table}

\subsection{Transient natural convection heatsink}
\label{sec34}

\begin{figure}[t]
    \centering
    \includegraphics[width=0.4\columnwidth]{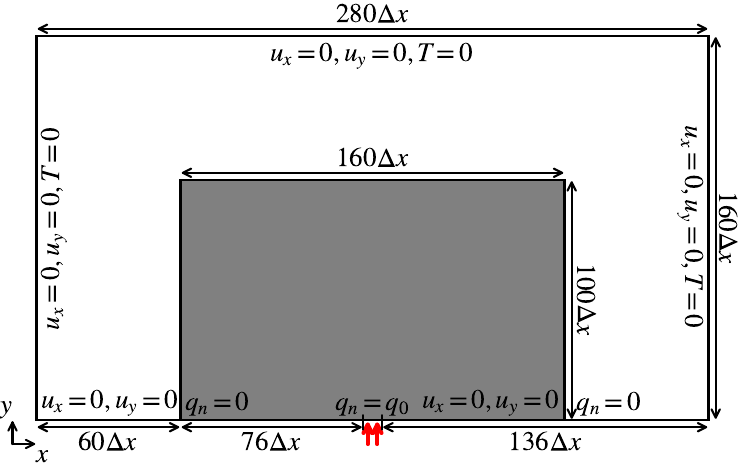}
    \caption{Design setting of transient natural convection heatsink}
    \label{fig:design_setting_heatsink}
\end{figure}

As the final example, we present a heat sink design under transient conditions to illustrate an unsteady natural convection problem.
Fig.~\ref{fig:design_setting_heatsink} shows the dimensions and boundary conditions.
The objective of this problem is to minimize the temperature on the heated boundary, as expressed in Eq.~\eqref{eq:objective_heatsink}.
The temperature is prescribed on the right, left and top walls, while the bottom wall is adiabatic, and the heated boundary is located at the center bottom, as shown in Fig.~\ref{fig:design_setting_heatsink}. 
\begin{equation}
    J=\frac{\int_\mathcal{I}\int_{\partial\mathcal{O}_\text{q}}Td\Gamma dt}{\int_\mathcal{I}\int_{\partial\mathcal{O}_\text{q}}d\Gamma dt}, \label{eq:objective_heatsink}
\end{equation}
where $\partial\mathcal{O}_\text{q}$ means the heated boundary.

In this problem, the external force based on the Boussinesq approximation is applied for the Navier-Stokes equation (Eq.~\eqref{eq:ns}) and the force is described as follows:
\begin{equation}
    G_\alpha=g_\alpha\beta\left(T-T_\text{ref}\right),\label{eq:boussinesqforce}
\end{equation}
where $g_\alpha$, $\beta$ and $T_\text{ref}$ are the gravitational acceleration vector, the volumetric expansion coefficient and the reference temperature, respectively.
Additionally, the solid area constraint, where $\gamma$ and $\mathcal{O}$ in Eq.~\eqref{eq:constraint_pipebend} are replaced by $1-\gamma$ and $D$, respectively, is applied, with its maximum value $V_\text{max}$ set to 0.5, in other words, 50\% of the area of the design domain.
In the heat sink design problem, it is widely known that numerous too fine elements are generated for the heat conduction dominant case.
To avoid generating these elements we use the Heaviside Projection Filter~\cite{wang2011projection}.
The filter radius $R$ is set to $2.4\Delta x$, and the steepness parameter $\beta$ is at first set to $1$ and is doubled per $100$ steps of the optimization loop or at the design variables fluctuation convergence.

In this numerical example, the Rayleigh number $Ra$ is defined by Eq.~\eqref{eq:rayleigh_heatsink}, and its value is set to $2\times10^5$,
\begin{equation}
    Ra=\frac{g\beta\Delta TL^3}{\nu K_\text{f}}. \label{eq:rayleigh_heatsink}
\end{equation}
Then, $g$, $\beta$, $\Delta T$ and $L$ are the gravitational acceleration, the volumetric expansion coefficient, the reference temperature difference and the characteristic length, and $\Delta T$ and $L$ are set to $1$ and $160\Delta x$, respectively.
Besides, the Prandtl number $Pr$ is set to $6$.
The thermal diffusivity of the solid $K_s$ is $10$ times larger than that of the fluid $K_f$.
The amount of the heat flux $q_0$ applied on the heated boundary is set to $1.0\times10^{-2}$.

\begin{figure}[t]
    \centering
    \begin{minipage}[b]{0.33\columnwidth}
        \centering
        \includegraphics[width=0.9\columnwidth]{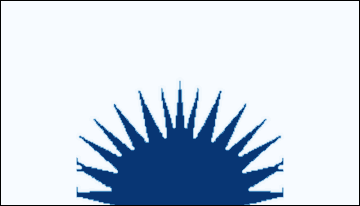}
        \subcaption{$t_1=5\times10^4$}
    \end{minipage}
    \begin{minipage}[b]{0.33\columnwidth}
        \centering
        \includegraphics[width=0.9\columnwidth]{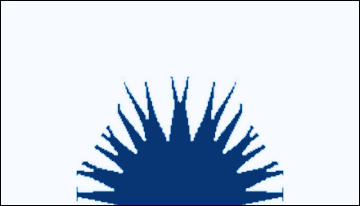}
        \subcaption{$t_1=7.5\times10^4$}
    \end{minipage}
    \begin{minipage}[b]{0.33\columnwidth}
        \centering
        \includegraphics[width=0.9\columnwidth]{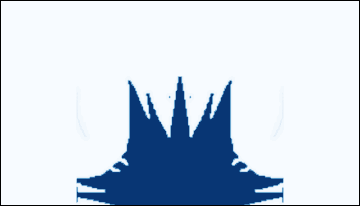}
        \subcaption{$t_1=1\times10^5$}
    \end{minipage}
    \caption{Optimized configuration of transient natural convection heatsink for each final time $t_1$}
    \label{fig:optimized_shape_heatsink}
\end{figure}

\begin{figure}[t]
    \centering
    \begin{minipage}[t]{0.07\columnwidth}
        \centering
        \includegraphics[width=1.0\columnwidth]{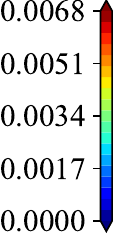}
    \end{minipage}
    \begin{minipage}[t]{0.28\columnwidth}
        \centering
        \includegraphics[width=1.0\columnwidth]{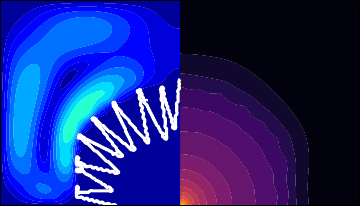}
        \subcaption{$t_1=5\times10^4$}
    \end{minipage}
    \begin{minipage}[t]{0.28\columnwidth}
        \centering
        \includegraphics[width=1.0\columnwidth]{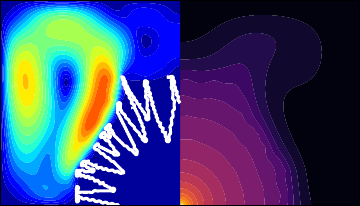}
        \subcaption{$t_1=7.5\times10^4$}
    \end{minipage}
    \begin{minipage}[t]{0.28\columnwidth}
        \centering
        \includegraphics[width=1.0\columnwidth]{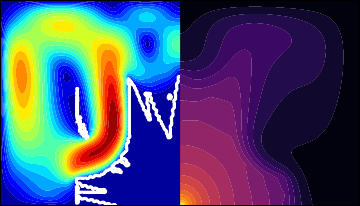}
        \subcaption{$t_1=1\times10^5$}
    \end{minipage}
    \begin{minipage}[t]{0.05\columnwidth}
        \centering
        \includegraphics[width=1.0\columnwidth]{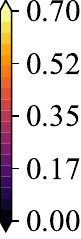}
    \end{minipage}
    \caption{Velocity and temperature distributions on optimized configuration for each final time $t_1$}
    \label{fig:velocity_temperature_heatsink}
\end{figure}

Fig.~\ref{fig:optimized_shape_heatsink} show the optimized configurations for the final time $t_1$ being $5\times10^4$, $7.5\times10^4$ and $10^5$, respectively.
Additionally, the each left half part of Fig.~\ref{fig:velocity_temperature_heatsink} shows the flow velocity distribution for $t_1=5\times10^4$, $t_1=7.5\times10^4$ and $t_1=10^5$ against the each right half part shows the temperature distribution for each final time $t_1$, respectively. 

For shorter $t_1$, the optimized configuration features multiple needle-like elements, whereas larger $t_1$ results in smoother surfaces along the flow.
The optimized configurations under the natural convection depend on the balance of the effects from the conduction and the convection, and at the early time for the simulation the flow is not much developed, so for the short $t_1$ condition the effect from the conduction should be relatively strong in the sensitivity.
The similar trend regarding the effect of the final time $t_1$ has been observed in previous research~\cite{Tanabe2023}.

The memory required by ALKS, given an analysis domain of $141 \times 161$ grid points and $50000$ time steps using double-precision floating-point format, is shown in Table~\ref{tab:heatsink_memory_usage_ALKS}, whereas the memory required by ALBM under the same conditions is shown in Table~\ref{tab:heatsink_memory_usage_ALBM}. 
Consequently, for sensitivity analysis, ALKS achieves approximately 75\% less memory usage than ALBM.
\begin{table}[t]
    \centering
    \begin{minipage}{\textwidth}
        \centering
        \caption{Memory requirements for transient natural convection heatsink using ALKS}
        \begin{tabular}{|p{0.45\columnwidth}|c|c|} \hline
            Variable name                                                                                                                                                                                                                                                                                                                                                                                & Array length              & Total amount                                       \\ \hline
            $\partial u_x/\partial x$,$\partial u_y/\partial x$,$\partial u_x/\partial y$,$\partial u_y/\partial y$, \newline $\partial\tilde{s}_{xx}/\partial x$,$\partial\tilde{s}_{xy}\partial x$,$\partial\tilde{s}_{yx}/\partial y$,$\partial\tilde{s}_{yy}/\partial y$, \newline $\partial T/\partial x$,$\partial T/\partial y$,$\partial\tilde{q}_x/\partial x$,$\partial\tilde{q}_y/\partial y$ & $n_x\times n_y$           & $141\times 161\times 12\times 8=3.88$MB            \\ \hline
            
            $\rho$,$q_x$,$q_y$,$\tilde{\rho}$,$\tilde{u}_x$,$\tilde{u}_y$,$\tilde{s}_{xx}$,$\tilde{s}_{xy}$,$\tilde{s}_{yx}$,$\tilde{s}_{yy}$,$\tilde{T}$,$\tilde{q}_x$,$\tilde{q}_y$                                                                                                                                                                                                                    & $n_x\times n_y\times 2$   & $141\times 161\times 2\times 13\times 8=4.72$MB    \\ \hline
            $u_x$,$u_y$,$T$                                                                                                                                                                                                                                                                                                                                                                              & $n_x\times n_y\times n_t$ & $141\times 161\times 50000\times 3\times 8=27.2$GB \\ \hline
        \end{tabular}
        \label{tab:heatsink_memory_usage_ALKS}
    \end{minipage}
    \begin{minipage}{\textwidth}
        \centering
        \caption{Memory requirements for transient natural convection heatsink using ALBM}
        \begin{tabular}{|c|c|c|} \hline
            Variable name                                                                                                                                    & Array length              & Total amount                                       \\ \hline
            $\rho$,$\tilde{\rho}$,$\tilde{m}_{x}$,$\tilde{m}_{y}$,$\tilde{T}$                                                                                & $n_x\times n_y$           & $141\times 161\times 5\times 8=908$kB              \\ \hline
            $q_x$,$q_y$,$\tilde{u}_x$,$\tilde{u}_y$,$\tilde{q}_{x}$,$\tilde{q}_y$,$f_0\dots f_8$,$\tilde{f}_0\dots\tilde{f}_8$,$\tilde{g}_0\dots\tilde{g}_8$ & $n_x\times n_y\times 2$   & $141\times 161\times 2\times 33\times 8=12.0$MB    \\ \hline
            
            $u_x$,$u_y$,$T$,$g_0\dots g_8$                                                                                                                   & $n_x\times n_y\times n_t$ & $141\times 161\times 50000\times 12\times 8=109$GB \\ \hline
        \end{tabular}
        \label{tab:heatsink_memory_usage_ALBM}
    \end{minipage}
\end{table}

\section{Conclusion}
\label{sec4}

In this study, we propose the \textit{Adjoint Lattice Kinetic Scheme} (ALKS) for topology optimization in fluid flow. 
In the proposed method, state fields are solved using the LKS, and sensitivities are calculated via the adjoint variable method, with adjoint equations derived from the discrete velocity Boltzmann equations and discretized and solved similarly to the LKS.
Since the boundary condition for the adjoint equations is difficult to solve directly due to inconsistency between the number of unknown variables and the number of equations, the method approximates the values employing the generalized inverse.
The sensitivity computed by ALKS is compared to that obtained by the finite difference approximation and, although small errors appear near the outer boundary or at the fluid-solid interface, the sensitivities are largely consistent and practically applicable.
The method is applied to the non-thermal and thermal steady-state and unsteady fluid problems, and in each case the results are verified by considering the physics and comparing them to previous studies for dependencies on parameters such as the Reynolds number or the final time of transient phenomena.
Furthermore, the memory usage in the sensitivity analysis with the proposed method is reduced by up to 75\% compared to the LBM for thermal unsteady fluid problems.
In this paper, numerical examples focus solely on 2D cases to quantitatively assess the computational efficiency of the proposed method using well-established problems.
However, the method can also be extended to 3D problems, where the conventional LBM approach would require significantly more fictitious particles. 
This extension to 3D presents a promising avenue for future work, potentially achieving an even greater reduction in memory usage with the proposed method.
 
\section{Replication of results}
The necessary information for replication of the results is presented in this paper.
The interested reader may contact the corresponding author for further implementation details.

\section*{Acknowledgements}
This work was supported by JSPS KAKENHI (GrantNo. 23K26018).

\section*{Appendix A: Derivation of the ALKS}
\label{appendixA}

We explain the derivation of the ALKS.
The sensitivity of the objective functional $J$, as expressed in Eq.~\eqref{eq:fj}, is calculated using the adjoint variable method.
The lattice Boltzman equations, shown in Eqs.~\eqref{eq:lbef}--\eqref{eq:lbeg_ex}, are derived from the discrete velocity Boltzman equations presented in Eqs.~\eqref{eq:boltzmanf} and \eqref{eq:boltzmang}:
\begin{align}
     & \text{Sh}\frac{\partial f_i}{\partial t}+c_{i\alpha}\frac{\partial f_i}{\partial x_\alpha}=-\frac{1}{\varepsilon_\text{f}}\left(f_i-f_i^\text{eq}\right)-3\Delta xw_ic_{i\alpha}\alpha_\gamma u_\alpha+3\Delta xw_ic_{i\alpha}G_\alpha,\label{eq:boltzmanf} \\
     & \text{Sh}\frac{\partial g_i}{\partial t}+c_{i\alpha}\frac{\partial g_i}{\partial x_\alpha}=-\frac{1}{\varepsilon_\text{g}}\left(g_i-g_i^\text{eq}\right)+\Delta xw_iQ.\label{eq:boltzmang}
\end{align}
Here, $\text{Sh}$ represents the Strouhal number.
$\varepsilon_\text{f}$ and $\varepsilon_\text{g}$ are dimensionless parameters of the same order as the Knudsen number, related to the dimensionless relaxation times $\tau_\text{f}$ and $\tau_\text{g}$, where $\tau_\text{f}=\varepsilon_\text{f}/\Delta x$ and $\tau_\text{g}=\varepsilon_\text{g}/\Delta x$, respectively.

The summation of the products of Eq.~\eqref{eq:boltzmanf} with the adjoint variable $\tilde{f}_i$ and Eq.~\eqref{eq:boltzmang} with $\tilde{g}_i$ is expressed as follows: 
\begin{align}
    R= & \int_{\mathcal{I}}\int_\mathcal{O}\sum_{i=0}^{Q-1}\tilde{f}_i\left\{\text{Sh}\frac{\partial f_i}{\partial t}+c_{i\alpha}\frac{\partial f_i}{\partial x_\alpha}+\frac{1}{\varepsilon_\text{f}}\left(f_i-f_i^{eq}\right)+3\Delta xw_ic_{i\alpha}\alpha_\gamma u_\alpha-3\Delta xw_ic_{i\alpha}G_\alpha\right\}d\Omega dt \label{eq:Lagrangian} \notag \\
       & +\int_{\mathcal{I}}\int_\mathcal{O}\sum_{i=0}^{Q-1}\tilde{g}_i\left\{\text{Sh}\frac{\partial g_i}{\partial t}+c_{i\alpha}\frac{\partial g_i}{\partial x_\alpha}+\frac{1}{\varepsilon_\text{g}}\left(g_i-g_i^{eq}\right)-\Delta xw_iQ\right\}d\Omega dt,
\end{align}
and the products of the boundary conditions and the adjoint variables $\tilde{\mu}_\alpha$, $\tilde{\lambda}$, $\tilde{\eta}$ and $\tilde{\kappa}$ expressed as follows:
\begin{align}
    R_\text{bu}= & \int_{\mathcal{I}}\int_{\partial\mathcal{O}_u}\tilde{\mu}_\alpha\left(u_\alpha-\bar{u}_\alpha\right)d\Gamma dt,                                                                  \\
    R_\text{bp}= & \int_{\mathcal{I}}\int_{\partial\mathcal{O}_p}\left\{\tilde{\lambda}\left(\rho-\bar{\rho}\right)+\tilde{\mu}_\text{s}\left(u_\text{s}-\bar{u}_\text{s}\right)\right\}d\Gamma dt, \\
    R_\text{bT}= & \int_{\mathcal{I}}\int_{\partial\mathcal{O}_T}\tilde{\eta}\left(T-\bar{T}\right)d\Gamma dt,                                                                                      \\
    R_\text{bq}= & \int_{\mathcal{I}}\int_{\partial\mathcal{O}_q}\tilde{\kappa}\left(q_\alpha n_\alpha+\bar{q}_\text{w}\right)d\Gamma dt
\end{align}
are summed and the Lagrangian $L$ in Eq.~\eqref{eq:lagrangian} is obtained.
\begin{equation}
    L=J+R+R_\text{bu}+R_\text{bp}+R_\text{bT}+R_\text{bq}\label{eq:lagrangian}
\end{equation}

Each term of Eq.~\eqref{eq:lagrangian} are differentiated by the design variable $\gamma$ and they become as follows:
\begin{align}
    \langle J^\prime,\delta\gamma\rangle=           & \int_{\mathcal{I}}\int_\mathcal{O}\left(\frac{\partial J_\mathcal{O}}{\partial\gamma}\delta\gamma+\sum_{i=0}^{Q-1}\frac{\partial J_\mathcal{O}}{\partial f_i}\delta f_i+\sum_{i=0}^{Q-1}\frac{\partial J_\mathcal{O}}{\partial g_i}\delta g_i\right)d\Omega dt \notag                                                                                                                                                                                                                                           \\                                                                                                                                                                                                                                          
                                                    & +\int_{\mathcal{I}}\int_{\partial\mathcal{O}}\left(\frac{\partial J_{\partial\mathcal{O}}}{\partial\gamma}\delta\gamma+\sum_{i=0}^{Q-1}\frac{\partial J_{\partial\mathcal{O}}}{\partial f_i}\delta f_i+\sum_{i=0}^{Q-1}\frac{\partial J_{\partial\mathcal{O}}}{\partial g_i}\delta g_i\right)d\Gamma dt,\label{eq:dfj}                                                                                                                                                                                          \\
    \langle R^\prime,\delta\gamma\rangle=           & \int_\mathcal{O}\sum_{i=0}^{Q-1}\text{Sh}\left[\tilde{f}_i\delta f_i\right]_{t_0}^{t_1}d\Omega+\int_{\mathcal{I}}\int_{\partial\mathcal{O}}\sum_{i=0}^{Q-1}n_\alpha\left\{\tilde{f}_i\delta_{\alpha\beta}-\frac{1}{\varepsilon_\text{f}}\Delta xA\left(\tilde{s}_{\alpha\beta}+\tilde{s}_{\beta\alpha}\right)\right\}c_{i\beta}\delta f_id\Gamma dt \notag                                                                                                                                                      \\
                                                    & +\int_{\mathcal{I}}\int_\mathcal{O}\sum_{i=0}^{Q-1}\left\{-\text{Sh}\frac{\partial\tilde{f}_i}{\partial t}-c_{i\alpha}\frac{\partial \tilde{f}_i}{\partial x_\alpha}+\frac{1}{\varepsilon_\text{f}}\left(\tilde{f}_i-\tilde{f}_i^{eq}\right)+3\Delta xc_{i\alpha}\alpha_\gamma \tilde{u}_\alpha-\frac{3}{\varepsilon_\text{g}}T\tilde{q}_\alpha c_{i\alpha}-3\Delta x\frac{\partial G_\alpha}{\partial f_i}\tilde{u}_\alpha-\Delta x\frac{\partial Q}{\partial f_i}\tilde{T}\right\}\delta f_id\Omega dt \notag \\
                                                    & +\int_{\mathcal{I}}\int_\mathcal{O}\left(3\Delta x\frac{\partial\alpha_\gamma}{\partial\gamma}u_\alpha\tilde{u}_\alpha-3\Delta x\frac{\partial G_\alpha}{\partial\gamma}\tilde{u}_\alpha\right)\delta\gamma d\Omega dt \notag                                                                                                                                                                                                                                                                                   \\
                                                    & +\int_\mathcal{O}\sum_{i=0}^{Q-1}\text{Sh}\left[\tilde{g}_i\delta g_i\right]_{t_0}^{t_1}d\Omega+\int_{\mathcal{I}}\int_{\partial\mathcal{O}}\sum_{i=0}^{Q-1}n_\alpha\left(\tilde{g}_i c_{i\alpha}-\frac{1}{\varepsilon_\text{g}}\Delta xB_\gamma\tilde{q}_\alpha\right)\delta g_id\Gamma dt \notag                                                                                                                                                                                                              \\
                                                    & +\int_{\mathcal{I}}\int_\mathcal{O}\sum_{i=0}^{Q-1}\left\{-\text{Sh}\frac{\partial\tilde{g}_i}{\partial t}-c_{i\alpha}\frac{\partial \tilde{g}_i}{\partial x_\alpha}+\frac{1}{\varepsilon_\text{g}}\left(\tilde{g}_i-\tilde{g}_i^{eq}\right)-3\Delta x\frac{\partial G_\alpha}{\partial g_i}\tilde{u}_\alpha-\Delta x\frac{\partial Q}{\partial g_i}\tilde{T}\right\}\delta g_id\Omega dt \notag                                                                                                                \\
                                                    & +\int_{\mathcal{I}}\int_\mathcal{O}\left(-\frac{1}{\varepsilon_\text{g}}\Delta x\frac{\partial B_\gamma}{\partial\gamma}\frac{\partial T}{\partial x_\alpha}\tilde{q}_\alpha-\Delta x\frac{\partial Q}{\partial\gamma}\tilde{T}\right)\delta\gamma d\Omega dt,                                                                                                                                                                                                                                                  \\
    \langle R_\text{bu}^\prime,\delta\gamma\rangle= & \int_{\mathcal{I}}\int_{\partial\mathcal{O}_u}\sum_{i=0}^{Q-1}\tilde{\mu}_\alpha c_{i\alpha}\delta f_id\Gamma dt,                                                                                                                                                                                                                                                                                                                                                                                               \\
    \langle R_\text{bp}^\prime,\delta\gamma\rangle= & \int_{\mathcal{I}}\int_{\partial\mathcal{O}_p}\sum_{i=0}^{Q-1}\left(\tilde{\lambda}+\tilde{\mu}_\text{s} c_{i\text{s}}\right)\delta f_id\Gamma dt,                                                                                                                                                                                                                                                                                                                                                              \\
    \langle R_\text{bT}^\prime,\delta\gamma\rangle= & \int_{\mathcal{I}}\int_{\partial\mathcal{O}_T}\sum_{i=0}^{Q-1}\tilde{\eta}\delta g_id\Gamma dt,                                                                                                                                                                                                                                                                                                                                                                                                                 \\
    \langle R_\text{bq}^\prime,\delta\gamma\rangle= & \int_{\mathcal{I}}\int_{\partial\mathcal{O}_q}\tilde{\kappa}\left\{-\sum_{i=0}^{Q-1}Tc_{i\alpha}n_\alpha\delta f_i+\sum_{i=0}^{Q-1}\left(c_{i\alpha}-u_\alpha\right)n_\alpha\delta g_i\right\}d\Gamma dt.
\end{align}
Here, $\tilde{f}_i^{eq}$ and $\tilde{g}_i^{eq}$ are expressed in Eqs.~\eqref{eq:afeq} and \eqref{eq:ageq}, respectively.
Besides,
\begin{align}
     & \tilde{\rho}=\sum_{i=0}^{Q-1}w_i\tilde{f}_i,\label{eq:afp}                                 \\
     & \tilde{u}_\alpha=\sum_{i=0}^{Q-1}w_ic_{i\alpha}\tilde{f}_i,\label{eq:afu}                  \\
     & \tilde{s}_{\alpha\beta}=\sum_{i=0}^{Q-1}w_ic_{i\alpha}c_{i\beta}\tilde{f}_i,\label{eq:afs} \\
     & \tilde{T}=\sum_{i=0}^{Q-1}w_i\tilde{g}_i,\label{eq:agT}                                    \\
     & \tilde{q}_\alpha=\sum_{i=0}^{Q-1}w_ic_{i\alpha}\tilde{g}_i,\label{eq:agq}
\end{align}
and they are corresponding the macroscopic values of the LKS.
The differential of the Lagrangian $L$ is arranged so that the coefficients of $\delta f_i$ and $\delta g_i$ become zero, yielding the following adjoint equations:
\begin{align}
     & -\text{Sh}\frac{\partial\tilde{f}_i}{\partial t}-c_{i\alpha}\frac{\partial \tilde{f}_i}{\partial x_\alpha}=-\frac{1}{\varepsilon_\text{f}}\left(\tilde{f}_i-\tilde{f}_i^{eq}\right)-3\Delta xc_{i\alpha}\alpha_\gamma \tilde{u}_\alpha+\frac{3}{\varepsilon_\text{g}}T\tilde{q}_\alpha c_{i\alpha}+3\Delta x\frac{\partial G_\alpha}{\partial f_i}\tilde{u}_\alpha+\Delta x\frac{\partial Q}{\partial f_i}\tilde{T}-\frac{\partial J_\mathcal{O}}{\partial f_i},\label{eq:aboltzmanf} \\
     & -\text{Sh}\frac{\partial\tilde{g}_i}{\partial t}-c_{i\alpha}\frac{\partial \tilde{g}_i}{\partial x_\alpha}=-\frac{1}{\varepsilon_\text{g}}\left(\tilde{g}_i-\tilde{g}_i^{eq}\right)+3\Delta x\frac{\partial G_\alpha}{\partial g_i}\tilde{u}_\alpha+\Delta x\frac{\partial Q}{\partial g_i}\tilde{T}-\frac{\partial J_\mathcal{O}}{\partial g_i},\label{eq:aboltzmang}
\end{align}
\begin{align}
     & \tilde{f}_i|_{t=t_1}=0,\label{eq:afinitial}                                                                                                                                                                                                                                                                                                                                                                          \\
     & \tilde{g}_i|_{t=t_1}=0,\label{eq:aginitial}                                                                                                                                                                                                                                                                                                                                                                          \\
     & n_\alpha\left\{\tilde{f}_i\delta_{\alpha\beta}-\frac{1}{\varepsilon_\text{f}}\Delta xA\left(\tilde{s}_{\alpha\beta}+\tilde{s}_{\beta\alpha}\right)\right\}c_{i\beta}+\tilde{\mu}_\alpha c_{i\alpha}+\frac{\partial J_{\partial\mathcal{O}}}{\partial f_i}=0                                                       & \left(\text{on }\partial\mathcal{O}_u\bigcap\partial\mathcal{O}_T\right),                        \\
     & n_\alpha\left\{\tilde{f}_i\delta_{\alpha\beta}-\frac{1}{\varepsilon_\text{f}}\Delta xA\left(\tilde{s}_{\alpha\beta}+\tilde{s}_{\beta\alpha}\right)\right\}c_{i\beta}+\tilde{\lambda}+\tilde{\mu}_\text{s}c_{i\text{s}}+\frac{\partial J_{\partial\mathcal{O}}}{\partial f_i}=0                                    & \left(\text{on }\partial\mathcal{O}_p\bigcap\partial\mathcal{O}_T\right),\label{eq:afboundaryp}  \\
     & n_\alpha\left(\tilde{g}_i c_{i\alpha}-\frac{1}{\varepsilon_\text{g}}\Delta xB_\gamma\tilde{q}_\alpha\right)+\tilde{\eta}+\frac{\partial J_{\partial\mathcal{O}}}{\partial g_i}=0                                                                                                                                  & \left(\text{on }\partial\mathcal{O}_T\right),\label{eq:agboundaryT}                              \\
     & n_\alpha\left\{\tilde{f}_i\delta_{\alpha\beta}-\frac{1}{\varepsilon_\text{f}}\Delta xA\left(\tilde{s}_{\alpha\beta}+\tilde{s}_{\beta\alpha}\right)\right\}c_{i\beta}+\tilde{\mu}_\alpha c_{i\alpha}-\tilde{\kappa}Tc_{i\alpha}n_\alpha+\frac{\partial J_{\partial\mathcal{O}}}{\partial f_i}=0                    & \left(\text{on }\partial\mathcal{O}_u\bigcap\partial\mathcal{O}_q\right),\label{eq:afboundaryuq} \\
     & n_\alpha\left\{\tilde{f}_i\delta_{\alpha\beta}-\frac{1}{\varepsilon_\text{f}}\Delta xA\left(\tilde{s}_{\alpha\beta}+\tilde{s}_{\beta\alpha}\right)\right\}c_{i\beta}+\tilde{\lambda}+\tilde{\mu}_\text{s}c_{i\text{s}}-\tilde{\kappa}Tc_{i\alpha}n_\alpha+\frac{\partial J_{\partial\mathcal{O}}}{\partial f_i}=0 & \left(\text{on }\partial\mathcal{O}_p\bigcap\partial\mathcal{O}_q\right),\label{eq:afboundarypq} \\
     & n_\alpha\left(\tilde{g}_i c_{i\alpha}-\frac{1}{\varepsilon_\text{g}}\Delta xB_\gamma\tilde{q}_\alpha\right)+\tilde{\kappa}\left(c_{i\alpha}-u_\alpha\right)n_\alpha+\frac{\partial J_{\partial\mathcal{O}}}{\partial g_i}=0                                                                                       & \left(\text{on }\partial\mathcal{O}_q\right).\label{eq:agboundaryq}
\end{align}
Eqs.~\eqref{eq:aboltzmanf} and \eqref{eq:aboltzmang} are discretized in the same manner as the LKS, and are expressed as follows:
\begin{align}
    \tilde{f}_i^*\left(\bm{x},t-\Delta t\right)= & \tilde{f}_i^\text{eq}\left(\bm{x}+\Delta x\bm{c}_i,t\right),\label{eq:albef}                                                                                                                                                                                                                                                   \\
    \tilde{f}_i\left(\bm{x},t-\Delta t\right)=   & \tilde{f}_i^*\left(\bm{x},t-\Delta t\right)-3\Delta xc_{i\alpha}\alpha_\gamma\left(\bm{x}\right)\tilde{u}_\alpha\left(\bm{x},t-\Delta t\right)+3T\left(\bm{x},t-\Delta t\right)\tilde{q}_\alpha\left(\bm{x},t-\Delta t\right)c_{i\alpha} \notag                                                                                \\
                                                 & +3\Delta x\frac{\partial G_\alpha}{\partial f_i}\left(\bm{x},t-\Delta t\right)\tilde{u}_\alpha\left(\bm{x},t\right)+\Delta x\frac{\partial Q}{\partial f_i}\left(\bm{x},t-\Delta t\right)\tilde{T}\left(\bm{x},t-\Delta t\right)-\frac{\partial J_\mathcal{O}}{\partial f_i}\left(\bm{x},t-\Delta t\right),\label{eq:albef_ex} \\
    \tilde{g}_i^*\left(\bm{x},t-\Delta t\right)= & \tilde{g}_i^\text{eq}\left(\bm{x}+\Delta x\bm{c}_i,t\right),\label{eq:albeg}                                                                                                                                                                                                                                                   \\
    \tilde{g}_i\left(\bm{x},t-\Delta t\right)=   & \tilde{g}_i^*\left(\bm{x},t-\Delta t\right)+3\Delta x\frac{\partial G_\alpha}{\partial g_i}\left(\bm{x},t-\Delta t\right)\tilde{u}_\alpha\left(\bm{x},t-\Delta t\right) \notag                                                                                                                                                 \\
                                                 & +\Delta x\frac{\partial Q}{\partial g_i}\left(\bm{x},t-\Delta t\right)\tilde{T}\left(\bm{x},t-\Delta t\right)-\frac{\partial J_\mathcal{O}}{\partial g_i}\left(\bm{x},t-\Delta t\right).\label{eq:albeg_ex}
\end{align}
Thus, the equations updating the macroscopic values of the adjoint variables are provided as follows:
\begin{align}
    \tilde{\rho}\left(\bm{x},t\right)=            & \sum_{i=0}^{Q-1}w_i\tilde{f}_i^\text{eq}\left(\bm{x}+\bm{c}_i\Delta x,t+\Delta t\right),                      \\
    \tilde{u}_\alpha\left(\bm{x},t\right)=        & \sum_{i=0}^{Q-1}w_ic_{i\alpha}\tilde{f}_i^\text{eq}\left(\bm{x}+\bm{c}_i\Delta x,t+\Delta t\right),           \\
    \tilde{s}_{\alpha\beta}\left(\bm{x},t\right)= & \sum_{i=0}^{Q-1}w_ic_{i\alpha}c_{i\beta}\tilde{f}_i^\text{eq}\left(\bm{x}+\bm{c}_i\Delta x,t+\Delta t\right), \\
    \tilde{T}\left(\bm{x},t\right)=               & \sum_{i=0}^{Q-1}w_i\tilde{g}_i^\text{eq}\left(\bm{x}+\bm{c}_i\Delta x,t+\Delta t\right),                      \\
    \tilde{q}_\alpha\left(\bm{x},t\right)=        & \sum_{i=0}^{Q-1}w_ic_{i\alpha}\tilde{g}_i^\text{eq}\left(\bm{x}+\bm{c}_i\Delta x,t+\Delta t\right).
\end{align}

\section*{Conflict of interest}
The authors declare that they have no conflict of interest.

\bibliographystyle{spbasic}
\bibliography{reference}

\end{document}